\documentclass[a4paper,11pt]{amsart}

\usepackage{a4wide}
\usepackage{tikz}
\usepackage{pgfplots}

\usepackage[foot]{amsaddr}

\newcommand{\SlopeTriangle}[6]
{
    \pgfplotsextra
    {
        \pgfkeysgetvalue{/pgfplots/xmin}{\xmin}
        \pgfkeysgetvalue{/pgfplots/xmax}{\xmax}
        \pgfkeysgetvalue{/pgfplots/ymin}{\ymin}
        \pgfkeysgetvalue{/pgfplots/ymax}{\ymax}

        \pgfmathsetmacro{\xArel}{#1}
        \pgfmathsetmacro{\yArel}{#3}
        \pgfmathsetmacro{\xBrel}{#1-#2}
        \pgfmathsetmacro{\yBrel}{\yArel}
        \pgfmathsetmacro{\xCrel}{\xArel}
        \pgfmathsetmacro{\lnxB}{\xmin*(1-(#1-#2))+\xmax*(#1-#2)} % in [xmin,xmax].
        \pgfmathsetmacro{\lnxA}{\xmin*(1-#1)+\xmax*#1} % in [xmin,xmax].
        \pgfmathsetmacro{\lnyA}{\ymin*(1-#3)+\ymax*#3} % in [ymin,ymax].
        \pgfmathsetmacro{\lnyC}{\lnyA+#4*(\lnxA-\lnxB)}
        \pgfmathsetmacro{\yCrel}{\lnyC-\ymin)/(\ymax-\ymin)} % THE IMPROVED EXPRESSION WITHOUT 'DIMENSION TOO LARGE' ERROR.

        \coordinate (A) at (rel axis cs:\xArel,\yArel);
        \coordinate (B) at (rel axis cs:\xBrel,\yBrel);
        \coordinate (C) at (rel axis cs:\xCrel,\yCrel);

        \draw[#6]   (A)--
                    (B)--
                    (C)-- node[anchor=east] {#5}
                    cycle;
    }
}

\usepackage[breaklinks,bookmarks=false]{hyperref}
\hypersetup{%
colorlinks,%
linkcolor=blue,%
citecolor=blue,%
urlcolor=blue,%
plainpages=false,%
pdfwindowui=false,%
pdfstartview={FitH}%
}

\newtheorem{theorem}{Theorem}
\newtheorem{lemma}[theorem]{Lemma}

\newtheorem{corollary}[theorem]{Corollary}
\newtheorem{remark}[theorem]{Remark}

\numberwithin{equation}{section}
\numberwithin{theorem}{section}

\usepackage{amssymb}
\usepackage{mathrsfs}

\newcommand{\eq}{:=}

\newcommand{\grad}{\boldsymbol \nabla}
\renewcommand{\div}{\grad \cdot}
\newcommand{\curl}{\grad \times}

\newcommand{\ccurl}{\boldsymbol{\operatorname{curl}}}
\newcommand{\ddiv}{\operatorname{div}}

\newcommand{\jmp}[1]{\,[\![#1]\!]}

\newcommand{\BA}{\boldsymbol A}
\newcommand{\BB}{\boldsymbol B}

\newcommand{\BD}{\boldsymbol D}
\newcommand{\BE}{\boldsymbol E}

\newcommand{\BH}{\boldsymbol H}
\newcommand{\BI}{\boldsymbol I}
\newcommand{\BJ}{\boldsymbol J}

\newcommand{\BL}{\boldsymbol L}

\newcommand{\BV}{\boldsymbol V}
\newcommand{\BW}{\boldsymbol W}

\newcommand{\ba}{\boldsymbol a}
\newcommand{\bb}{\boldsymbol b}

\newcommand{\bn}{\boldsymbol n}
\newcommand{\bo}{\boldsymbol o}

\newcommand{\bv}{\boldsymbol v}
\newcommand{\bw}{\boldsymbol w}
\newcommand{\bx}{\boldsymbol x}
\newcommand{\by}{\boldsymbol y}

\newcommand{\CF}{\mathcal F}

\newcommand{\CM}{\mathcal M}
\newcommand{\CN}{\mathcal N}

\newcommand{\CP}{\mathcal P}

\newcommand{\CT}{\mathcal T}

\newcommand{\CV}{\mathcal V}

\newcommand{\LC}{\mathscr C}

\newcommand{\LP}{\mathscr P}

\newcommand{\BCN}{\boldsymbol{\CN}}

\newcommand{\BCP}{\boldsymbol{\CP}}

\newcommand{\RT}{\boldsymbol{RT}}
\newcommand{\N}{\boldsymbol{N}}

\newcommand{\ee}{\boldsymbol \varepsilon}
\newcommand{\mm}{\boldsymbol \mu}
\newcommand{\cc}{\boldsymbol \chi}

\newcommand{\bphi}{\boldsymbol \phi}
\newcommand{\bpsi}{\boldsymbol \psi}

\newcommand{\osc}{\operatorname{osc}}

\newcommand{\oma}{\omega^{\ba}}
\newcommand{\gma}{\Gamma^{\ba}}
\newcommand{\pa}{\psi^{\ba}}
\newcommand{\CTa}{\CT_h^{\ba}}

\newcommand{\GT}{\Gamma_{\rm T}}
\newcommand{\GN}{\Gamma_{\rm N}}

\newcommand{\LBA}{\boldsymbol \Lambda}

\newcommand{\ND}{\BCN}

\newcommand{\Cconta}{C_{{\rm cont},\ba}}
\newcommand{\Csta}{C_{{\rm st},\ba}}

\newcommand{\bzero}{\bo}

\newcommand{\bxi}{\boldsymbol \xi}

\newcommand{\CcontK}{C_{{\rm cont},K}}
\newcommand{\CstK}{C_{{\rm st},K}}
\newcommand{\tK}{\widetilde K}

\newcommand{\merrH}{\text{err}_{\rm H}}
\newcommand{\merrB}{\text{err}_{\rm B}}
\newcommand{\mest }{\eta}
\newcommand{\ndofs}{N_{\rm dofs}}

\newcommand{\errH}{$\merrH$}
\newcommand{\errB}{$\merrB$}
\newcommand{\est }{$\mest$}
\newcommand{\ndof}{$\ndofs$}

\title[A posteriori estimates for mixed discretizations of curl-curl problems]{An equilibrated estimator for mixed finite element discretizations of the curl-curl problem}
\author{T. Chaumont-Frelet$^\dagger$}

\address{\vspace{-.5cm}}
\address{\noindent \tiny \textup{$^\dagger$Inria Universit\'e C\^ote d'Azur, LJAD, CNRS}}

\begin{document}

\maketitle

\begin{abstract}
We propose a new a posteriori error estimator for mixed finite element discretizations
of the curl-curl problem. This estimator relies on a Prager--Synge inequality, and therefore
leads to fully guaranteed constant-free upper bounds on the error. The estimator is
also locally efficient and polynomial-degree-robust. The construction is based on patch-wise
divergence-constrained minimization problems, leading to a cheap embarrassingly parallel
algorithm. Crucially, the estimator operates without any assumption on the topology of
the domain, and unconventional arguments are required to establish the reliability estimate.
Numerical examples illustrate the key theoretical results, and suggest that the estimator is
suited for mesh adaptivity purposes.

\vspace{.5cm}
\noindent
{\sc Key words.}
a posteriori error estimate;
electromagnetics
finite element method;
high order method;
potential reconstruction;
Prager--Synge
\end{abstract}

\section{Introduction}

This work develops an equilibrated a posteriori error estimator for mixed finite
element discretizations of the curl-curl problem. The curl-curl equation is the
prototypical elliptic PDE in $\BH(\ccurl)$, and constitutes the basic model problem
for magnetostatics. In contrast to the $\BH(\ccurl)$ setting, the development of
equilibrated estimators in $H^1$ with application to, e.g., electrostatics, is much
more advanced. This introduction reviews the key concepts of equilibrated estimators
in $H^1$, highlights the challenges that arise in $\BH(\ccurl)$, presents the construction
of the new estimator and summarizes the key results of the present work.

\subsection{Equilibrated estimators in $H^1$}

A density of static charges $\rho: \Omega \to \mathbb R$ generates
an electric field $\BE: \Omega \to \mathbb R^3$. The electric
field is linked to the electric displacement $\BD = \ee\BE$, where
$\ee$ is the electric permittivity tensor, and Gauss' law
ensures that $\div \BD = \rho$. Assuming a steady state,
Faraday's law stipulates the existence an electric potential
$\phi: \Omega \to \mathbb R$ such that $\BE = \grad \phi$. This
leads to the electrostatic problem of finding $\phi$ such that
\begin{equation}
\label{eq_electro_strong}
\left \{
\begin{array}{rcll}
-\div (\ee\grad \phi) &=& \rho & \text{ in } \Omega,
\\
\phi &=& 0 & \text{ on } \GT,
\\
\grad \phi \cdot \bn &=& 0 & \text{ on } \GN,
\end{array}
\right .
\end{equation}
where the partition $\{\GT,\GN\}$ of the boundary $\partial \Omega$
depends on the properties of the materials surrounding $\Omega$.

If $\BE_h$ is any square-integrable vector field, then the
following generalized Prager--Synge inequality holds true
\begin{equation}
\label{eq_prager_synge_electro}
\|\BE-\BE_h\|_{\ee,\Omega}^2
\leq
\min_{\widetilde \phi \in H^1_{\GT}(\Omega)}
\|\BE_h-\grad \widetilde \phi\|_{\ee,\Omega}^2
+
\min_{\substack{\widetilde \BD \in \BH_{\GN}(\ddiv,\Omega) \\ \div \widetilde \BD = \rho}}
\|\BE_h-\ee^{-1} \widetilde \BD\|_{\ee,\Omega}^2,
\end{equation}
as shown for instance in \cite[Theorem 3.3]{ern_vohralik_2015a}. The name Prager--Synge
is after the seminal work \cite{prager_synge_1947a}, where Prager and Synge proved
\eqref{eq_prager_synge_electro} in the context of elastostatics for fields
$\BE_h = \grad \phi_h$, $\phi_h \in H^1_{\GT}(\Omega)$, hence with the first term vanishing.
It is noteworthy that the two terms of the right-hand of \eqref{eq_prager_synge_electro}
respectively quantify the inability of $\BE_h$ to satisfy Gauss' and Faraday's laws.

When numerically solving \eqref{eq_electro_strong} with a conforming finite element method,
one of the two terms in the right-hand side of \eqref{eq_prager_synge_electro} automatically
vanishes by construction. The remainding term can then be bounded using an error estimator,
leading to guaranteed error bounds \cite{ern_vohralik_2015a}.
%?

%?
Specifically, when Lagrange elements are used to approximate \eqref{eq_electro_strong}
in primal form, a conforming electric potential $\phi_h \in H^1_{\GT}(\Omega)$ is computed,
and the resulting electric field approximation is $\BE_h \eq \grad \phi_h$.
The first term in the right-hand side of \eqref{eq_prager_synge_electro} then vanishes, and
the second term is controlled using an equilibrated flux reconstruction whereby a field
$\BD_h$ such that $\div \BD_h = \rho$ is computed
\cite{
ainsworth_oden_2000a,%
destuynder_metivet_1999a,%
ern_vohralik_2015a,%
ladeveze_leguillon_1983a,%
luce_wohlmuth_2004a,%
nicaise_witowski_wohlmuth_2008a}.

On the other hand, when employing Raviart--Thomas elements for the mixed form
of \eqref{eq_electro_strong}, an electric displacement $\BD_h$ satisfying $\div \BD_h = \rho$
is immediately computed. In this case, the second term in the right-hand
side of \eqref{eq_prager_synge_electro} vanishes. The first term is then estimated by building
a field $\phi_h$ such that $\ee\BD_h-\grad \phi_h$ is small: this process is called a
potential reconstruction
\cite{ainsworth_2008a,ern_vohralik_2015a,ern_vohralik_2021a,vohralik_2007a}.

In practice, for both primal and mixed forms, the identity $\div \BD_h = \rho$ only holds
point-wise if $\rho$ is piecewise polynomial on the finite element mesh. In the general case,
an additional term corresponding to the approximation of $\rho$ is added to
\eqref{eq_prager_synge_electro}. This term is fully computable, and when $\rho$ is piecewise
smooth, it is actually of higher-order and is usually called a ``data oscillation'' term.

A chief question is then: How such equilibrated flux and potential should be constructed?
For primal discretizations, several approaches have been proposed in the past for the flux
reconstruction
\cite{
ainsworth_oden_2000a,%
ladeveze_leguillon_1983a,%
luce_wohlmuth_2004a,%
nicaise_witowski_wohlmuth_2008a},
and in this work, we will follow the technique introduced in \cite{destuynder_metivet_1999a}.
Specifically, instead of solving the global minimization problem
\begin{equation*}
\min_{\substack{\BD_h \in \BH_{\GN}(\ddiv,\Omega) \cap \RT_p(\CT_h) \\ \div \BD_h = \rho}}
\|\BE_h-\ee^{-1} \BD_h\|_{\ee,\Omega}^2,
\end{equation*}
with Raviart--Thomas elements, where $p$ is the polynomial degree of the Lagrange finite
elements, the localized version
\begin{equation}
\label{eq_intro_BDa}
\BD_h^{\ba}
\eq
\min_{\substack{\bv_h \in \BH_0(\ddiv,\oma) \cap \RT_{p+1}(\CTa)
\\
\div \bv_h = \pa \rho-\grad \pa \cdot \grad \phi_h}}
\|\BE_h-\ee^{-1} \bv_h\|_{\ee,\oma}^2,
\end{equation}
is considered for each vertex $\ba$ of the mesh, where $\pa$ is the associated
hat function (these notations are rigorously introduced in Section
\ref{section_vertex_patches} below). Crucially, the compatibility condition
implied by Stokes' formula in \eqref{eq_intro_BDa} is satisfied due to Galerkin
orthogonality. Indeed $\phi_h$ is the Galerkin finite element solution and $\pa$
is a lowest-order Lagrange finite element function. One easily checks that the
corresonding field
\begin{equation*}
\BD_h \eq \sum_{\ba \in \CV_h} \BD_h^{\ba}
\end{equation*}
is an equilibrated flux, so that \eqref{eq_prager_synge_electro} leads to a constant-free
error estimate.

Similar localization stratgies involving the hat functions of the finite element
mesh have also been introduced for the potential reconstruction, and we refer the
reader to \cite{ern_vohralik_2015a} for more details.

Another key aspect of equilibrated estimators built from local minimization problems
is that they are ``polynomial-degree-robust''
\cite{braess_pillwein_schoberl_2009a,ern_vohralik_2015a,ern_vohralik_2021a}.
It means that the estimator is also a lower bound to the error up to a generic
constant, and that this constant does not depend on the polynomial degree $p$ of
the finite element space. Therefore, equilibrated estimators are (at least theoretically)
more attractive than traditional residual-based estimators for high-order finite
element methods and/or $hp$-adaptive algorithms.

\subsection{Prager--Synge estimates in $\BH(\ccurl)$}

In this work, we are rather interested in magnetostatics \cite{griffiths_1999a}.
In this case, a magnetic field $\BH: \Omega \to \mathbb R^3$ is generated by a
static current $\BJ: \Omega \to \mathbb R^3$ with $\div \BJ = 0$. The magnetic
field is linked to the induction field $\BB$ through the constitutive relation $\BB = \mm\BH$,
where $\mm$ is the magnetic permeability tensor. Amp\`ere's law and Gauss'
law for magnetism then respectively demand that
$\curl \BH = \BJ$ and $\BB = \curl \BA$ for some vector-potential
$\BA: \Omega \to \mathbb R^3$ with $\div \BA = 0$. The curl-curl equation is then obtained
by reformulating the problem in terms of $\BA$ as follows:
\begin{equation}
\label{eq_magneto_strong}
\left \{
\begin{array}{rcll}
\curl \left (\mm^{-1}\curl\BA\right ) &=& \BJ & \text{ in } \Omega,
\\
\div \BA &=& 0 & \text{ in } \Omega
\\
(\curl \BA) \times \bn &=& \bo & \text{ on } \GT,
\\
\BA \times \bn &=& \bo & \text{ on } \GN.
\end{array}
\right .
\end{equation}
Depending on the topology of $\Omega$, $\GT$ and $\GN$ a finite number of
additional gauge conditions may be imposed to uniquely define $\BA$. This
is detailed in Section \ref{section_cohomology} below.

The magnetostatic problem in \eqref{eq_magneto_strong} has a structure
similar to the electrostatic problem. As a result, it also enjoys a
generalized Prager--Synge inequality. Specifically, if $\BH_h$ is any
square-integrable vector field, we have
\begin{equation}
\label{eq_prager_synge_magneto_intro}
\|\BH-\BH_h\|_{\mm,\Omega}^2
\leq
\min_{\widetilde \BA \in \BH_{\GN}(\ccurl,\Omega)}
\|\BH_h-\mm^{-1}\curl \widetilde \BA\|_{\mm,\Omega}^2
+
\min_{\substack{\widetilde \BH \in \BH_{\GT}(\ccurl,\Omega) \\ \curl \widetilde \BH = \BJ}}
\|\BH_h-\widetilde \BH\|_{\mm,\Omega}^2,
\end{equation}
where, as a counterpart to \eqref{eq_prager_synge_electro}, the two terms in the right-hand
side respectively measure the failure of Amp\`ere's and Gauss' laws.

The curl-curl equation \eqref{eq_magneto_strong} can be immediately discretized
in primal form. In this case, an approximation $\BA_h$ of $\BA$ is computed with
N\'ed\'elec elements, and the divergence constraint is imposed with (the gradient of)
Lagrange elements as Lagrange multipliers. As for primal discretizations of
\eqref{eq_electro_strong}, in this case, we set $\BH_h \eq \mm^{-1} \BA_h$, and
the first term in the right-hand side of \eqref{eq_prager_synge_magneto_intro}
vanishes. The error is then estimated by reconstructing a field $\widetilde \BH$
such that $\curl \widetilde \BH = \BJ$, up to data oscillations. These curl-constrained
equilibrated reconstructions turn out to be substantially more complex than in the $H^1$
setting, and general polynomial degrees $p$ have only been handled very recently
\cite{braess_schoberl_2008a,
chaumontfrelet_vohralik_2022a,
gedicke_geevers_perugia_2020a,
gedicke_geevers_perugia_schoberl_2021a}.

Here, we focus on the discretization of the mixed form of \eqref{eq_magneto_strong}
which is not currently covered in the literature. Namely, $\BH$ is immediately discretized
with N\'ed\'elec elements, so that $\curl \BH_h = \BJ$ and the second term in the right-hand
side of \eqref{eq_prager_synge_magneto_intro} vanishes up to data oscillation. A potential
reconstruction is thus required to control the first term.

\subsection{A novel error estimator}

To analyze the efficiency of estimators based on a potential reconstruction
in $H^1$, one of the key tools is a broken Poincar\'e inequality for scalar
functions with jumps of vanishing mean value, see e.g.
\cite[Lemma 3.13 and Theorem 3.17]{ern_vohralik_2015a}
and \cite[Corollary 4.1]{ern_vohralik_2021a} as well as
\cite{brenner_2003a,knobloch_2001a,vohralik_2005a}. However,
to the best of the author's knowledge, such an inequality is
not available in the $\BH(\ccurl)$ context.

In this work, we therefore follow an alternative strategy. Specifically, instead
of constructing a potential $\BA_h$ with N\'ed\'elec elements and then
take its curl in \eqref{eq_prager_synge_magneto_intro}, we immediately construct a Raviart--Thomas
function $\BB_h$ in the range of the curl operator. $\BB_h$ is
then used in place of $\curl \BA_h$ in \eqref{eq_prager_synge_magneto_intro}.

If $\Omega$ is homotopy equivalent to a ball, and either $\GT$ or $\GN = \emptyset$,
then it is in fact sufficient to build a divergence-free field $\BB_h \in \BH_{\GN}(\ddiv,\Omega)$.
This observation leads us to consider a divergence-constrained minimization problem of the form
\begin{equation*}
\min_{\substack{\widetilde \BB \in \BH_{\GN}(\ddiv,\Omega) \\ \div \widetilde \BB = 0}}
\|\BH_h-\widetilde \mm^{-1}\widetilde \BB\|_{\mm,\Omega}^2,
\end{equation*}
which is very similar to the second term in the right-hand side of \eqref{eq_prager_synge_electro}.
In fact, the localization technique introduced in \cite{destuynder_metivet_1999a} and further
analyzed in \cite{braess_pillwein_schoberl_2009a,ern_vohralik_2015a} can be accommodated,
leading to the localized divergence-constrained problems
\begin{subequations}
\label{eq_estimator_intro}
\begin{equation}
\label{eq_local_problem_intro}
\BB_h^{\ba}
\eq
\min_{\substack{
\bv_h \in \BH_0(\ddiv,\oma) \cap \RT_{p+2}(\CTa)
\\
\div \bv_h = \grad \pa \cdot \BH_h
}}
\|\pa \BH_h-\mm^{-1}\bv_h\|_{\mm,\oma}^2.
\end{equation}
These local contributions are then assembled
\begin{equation}
\label{eq_equilibrated_field_intro}
\BB_h \eq \sum_{\ba \in \CV_h} \BB_h^{\ba}
\end{equation}
leading to the estimator
\begin{equation}
\label{eq_eta_intro}
\eta_K \eq \|\BH_h-\mm^{-1}\BB_h\|_{\mm,K},
\end{equation}
for all $K \in \CT_h$.
It is clear that $\div \BB_h = 0$, so that $\BB_h$ is an admissible field
to plug in \eqref{eq_prager_synge_magneto} when the cohomology of $\Omega$ is trivial,
thus leading to guaranteed error bounds.
\end{subequations}

The interesting (and perhaps surprising) part of this work, is that this simple procedure
is still valid without any assumption on the topology of $\Omega$, $\GT$ or $\GN$. Specifically,
not only does the localization procedure produce a divergence-free field $\BB_h$, but it also
guarantees that this field is in the range of the curl operator. This result is obtained by
checking that the flux of $\BB_h$ vanishes through all closed surfaces. This is a striking
result, since it is a global property that is not actively enforced in the construction
based on local problems.

\subsection{Main results}

The key theoretical results of this work may be summarized as follows.
The local problems \eqref{eq_local_problem_intro} are well-posed, and
the vector field $\BB_h$ constructed following \eqref{eq_equilibrated_field_intro}
always sits in the range of the curl operator, without any assumption on the topology.
The corresponding error estimator in \eqref{eq_eta_intro} is reliable up to data
oscillations and locally efficient:
\begin{equation}
\label{eq_properties_intro}
\|\BH-\BH_h\|_{\mm,\Omega}^2
\leq
\sum_{K \in \CT_h} \eta_K^2 + \osc^2,
\qquad
\eta_K \lesssim
\|\BH-\BH_h\|_{\mm,\widetilde K},
\end{equation}
for all $K \in \CT_h$. In addition, the hidden constant in the lower bound does not depend on
the polynomial degree $p$, so that the estimator is polynomial-degree-robust. These estimates
are rigorously stated and established in Theorem \ref{theorem_upper_bound} and Corollary
\ref{corollary_lower_bound} below.

We also illustrate the theoretical findings with numerical examples. Crucially, we
observe that \eqref{eq_properties_intro} indeed sharply holds numerically. We also
investigate the accuracy of the reconstructed field $\BB_h$ as opposed to $\mm\BH_h$
to approximate the magnetic induction $\BB$. Our observation is that although $\BB_h$
does not exhibit super convergence properties, it is $2$ to $5$ times more accurate
than $\mm\BH_h$ in our examples. We finally employ the estimator for adaptive mesh
refinements, and observe optimal convergence rates.

\subsection{Outline}

The remainder of this work is organized as follows. Section \ref{section_continuous_setting}
properly states our model problem, whereas Section \ref{section_discrete_setting} introduces
the discretization setting. Sections \ref{section_reliability} and \ref{section_efficiency}
respectively contain the reliability and efficiency proofs. We present numerical examples
in Section \ref{section_numerical_examples}, before drawing concluding remarks in
Section \ref{section_conclusion}.

\section{Continuous setting}
\label{section_continuous_setting}

\subsection{Domain and coefficient}

We consider a polyhedral domain $\Omega \subset \mathbb R^3$ with a Lipschitz boundary
$\partial \Omega$. $\partial \Omega$ is split into two disjoint, relatively open, and polytopal
components $\GN$ and $\GT$ with Lipschitz boundaries (the case where these manifolds do not have
boundaries is allowed).

We also fix a symmetric matrix-valued coefficient $\mm: \Omega \to \mathbb R^{3 \times 3}$.
Classically, we require that $\mm$ is uniformly elliptic and bounded, and for the sake of
simplicity, $\cc \eq \mm^{-1}$ will denote the inverse of $\mm$. We further assume that
there exists a partition $\LP$ of $\Omega$ into non-overlapping polyhedral subsets $P$
such that for each $P \in \LP$, $\mm|_P$ and $\cc|_P$ are constant matrices.

The ``coefficient contrast'' will appear in the efficiency bounds of the proposed
estimator. For any open subset $\omega \subset \Omega$, it is defined as
\begin{equation*}
\LC_{\mm,\omega}
\eq
\underset{\bx \in \omega}{\text{ess sup}} \max_{\substack{\bxi \in \mathbb R^3 \\ |\bxi| = 1}}
\mm \bxi \cdot \bxi
\Big /
\underset{\bx \in \omega}{\text{ess inf}} \min_{\substack{\bxi \in \mathbb R^3 \\ |\bxi| = 1}}
\mm \bxi \cdot \bxi.
\end{equation*}

\subsection{Functional spaces}

For an open set $\omega \subset \mathbb R^3$, $L^2(\omega)$ is the usual Lebesgue space
of real-valued square integrable functions, and $\BL^2(\omega) \eq [L^2(\omega)]^3$ contains
vector-valued functions \cite{adams_fournier_2003a}. The standard inner-products
of both $L^2(\omega)$ and $\BL^2(\omega)$ are denote by $(\cdot,\cdot)_\omega$. We employ the
notation $\|{\cdot}\|_\omega$ for the usual norm of $\BL^2(\omega)$ associated with
$({\cdot},{\cdot})_\omega$. We will also frequently used the (equivalent) weighted norms
$\|{\cdot}\|_{\mm,\omega}^2 \eq (\mm{\cdot},{\cdot})_\omega$ and
$\|{\cdot}\|_{\cc,\omega}^2 \eq (\cc{\cdot},{\cdot})_\omega$.

For standard Sobolev spaces, we will use the notations
\begin{align*}
H^1(\omega)
&\eq
\left \{ v \in L^2(\omega); \; \grad v \in \BL^2(\omega) \right \},
\\
\BH(\ccurl,\omega)
&\eq
\left \{ \bv \in \BL^2(\omega); \; \curl \bv \in \BL^2(\omega) \right \},
\\
\BH(\ddiv,\omega)
&\eq
\left \{ \bv \in \BL^2(\omega); \; \div \bv \in L^2(\omega) \right \},
\end{align*}
where $\grad$, $\curl$ and $\div$ are the weak gradient, curl and divergence
operators defined in the sense of distributions, see
\cite{adams_fournier_2003a,girault_raviart_1986a}.

If $\gamma \subset \partial \omega$, we respectively denote by $H^1_\gamma(\omega)$,
$\BH_\gamma(\ccurl,\omega)$ and $\BH_\gamma(\ddiv,\omega)$ the closure of smooth
functions that vanish on $\gamma$ in $H^1(\omega)$, $\BH(\ccurl,\omega)$
and $\BH(\ddiv,\omega)$. Assuming that $\omega$ has a Lipschitz boundary and that $\gamma$
is relatively open, these spaces may be interpreted as containing functions with
vanishing trace, tangential trace and normal trace on $\gamma$, see \cite{fernandes_gilardi_1997a}.

Finally, we introduce the short-hand notation
$\LBA_\gamma(\omega) \eq \curl \BH_\gamma(\ccurl,\omega)$.

\subsection{Cohomology}
\label{section_cohomology}

The space $\LBA_{\GN}(\Omega)$ will play a key role in the following.
It is clear that for all $\bv \in \BH_{\GN}(\ccurl,\Omega)$,
$\bw \eq \curl \bv \in \BH_{\GN}(\ddiv,\Omega)$ with $\div \bw = 0$,
so that
\begin{equation*}
\LBA_{\GN}(\Omega)
\subset
\BH_{\GN}(\ddiv^0,\Omega)
\eq
\left \{
\bw \in \BH_{\GN}(\ddiv,\Omega)
\; | \;
\div \bw = 0
\right \}.
\end{equation*}
In fact, if $\Omega$ is homotopy equivalent to a ball and either $\GT = \emptyset$
or $\GN = \emptyset$, then the identity $\LBA_{\GN}(\Omega) = \BH_{\GN}(\ddiv^0,\Omega)$
holds true.

In the general case, a finite number of linear constraints must be
imposed on $\bw \in \BH_{\GN}(\ddiv^0,\Omega)$ to ensure the existence
of $\bv \in \BH_{\GN}(\ccurl,\Omega)$ such that $\bw = \curl \bv$.
Specifically, following \cite{assous_ciarlet_labrunie_2018a,duff_1952a,gross_kotiuga_2004a},
there exists a finite number of (relatively open) oriented surfaces
$\Sigma_j \subset \overline{\Omega}$ with unit normal vector $\bn_{\Sigma_j}$
and $\partial \Sigma_j \subset \overline{\GN}$, $1 \leq j \leq N$, such that
$\bw \in \LBA_{\GN}(\ddiv,\Omega)$ if and only if $\bw \in \BH_{\GN}(\ddiv^0,\Omega)$ and
\begin{equation}
\label{eq_condition_sigma}
\int_{\Sigma_j} \bw \cdot \bn_{\Sigma_j} = 0
\end{equation}
for $1 \leq j \leq N$. The integrals in the left-hand side of \eqref{eq_condition_sigma}
are called the ``periods'' of $\bw$. A function $\bw \in \BH_{\GN}(\ddiv^0,\Omega)$ thus
admits a vector potential if and only if its periods vanish.

\begin{remark}[Practical implementation]
The surfaces $\{\Sigma_j\}_{j=1}^N$ are only needed to develop
the theory. They are not used in the implementation and do not
explicitly enter the construction of the proposed error estimator.
\end{remark}

%%  \begin{remark}[Intuitive explanation]
%%  Let us first consider a smooth function $\bv \in \BC^1(\overline{\Omega})$
%%  with $\bv \times \bn = \bzero$ on $\GT$, and $\Sigma \subset \Omega$ is a
%%  smooth surface such that $\partial \Sigma \subset \overline{\GN}$, we have
%%  \begin{equation*}
%%  \int_{\Sigma} \curl \bv \cdot \bn_{\Sigma}
%%  =
%%  \int_{\partial \Sigma} \bv \cdot \bt_{\Sigma}
%%  =
%%  0.
%%  \end{equation*}
%%  Considering on the other hand $\bw \in \BC^0(\overline{\Omega})$
%%  with $\div \bw = 0$, we easily see that
%%  \begin{equation*}
%%  \int_{\partial \omega} \bw \cdot \bn
%%  =
%%  \int_{\omega} \div \bw
%%  =
%%  0
%%  \end{equation*}
%%  for all smooth subdomains $\omega \subset \Omega$, as a result
%%  we do have
%%  \begin{equation*}
%%  \int_{\Sigma} \bw \cdot \bn = 0
%%  \end{equation*}
%%  for all $\Sigma$ such that $\Sigma = \partial \omega$ for some subdomain $\omega \subset \Omega$.
%%  \end{remark}

\subsection{Model problem}

Given a right-hand side $\BJ \in \LBA_{\GT}(\Omega)$, our model problem is to find
$\BH \in \BH_{\GT}(\ccurl,\Omega)$ and $\BA \in \LBA_{\GT}(\Omega)$ such that
\begin{equation}
\label{eq_magneto_continuous}
\left \{
\begin{array}{rcl}
(\mm\BH,\bv)_\Omega + (\BA,\curl\bv)_\Omega &=& 0,
\\
(\curl \BH,\bw)_\Omega &=& (\BJ,\bw)_\Omega,
\end{array}
\right .
\end{equation}
for all $\bv \in \BH_{\GT}(\ccurl,\Omega)$ and $\bw \in \LBA_{\GT}(\Omega)$.
The system of equations in \eqref{eq_magneto_continuous} is a saddle-point
problem. By definition of $\LBA_{\GT}(\Omega)$, the so-called ``inf-sup''
condition is trivially satisfied, so that there exists a unique solution
\cite{brezzi_1974a}.

\section{Discrete setting}
\label{section_discrete_setting}

\subsection{Computational mesh}

Throughout this work, we consider a fixed mesh $\CT_h$ of $\Omega$
consisting of (open) tetrahedral elements $K$. The set of all mesh
vertices and faces are respectively denoted by $\CV_h$ and $\CF_h$,
and we further split $\CF_h$ into the set $\CF_h^{\rm e}$
of exterior faces $F \in \CF_h$ such that $F \subset \partial \Omega$
and the set interior faces $\CF_h^{\rm i} \eq \CF_h \setminus \CF_h^{\rm e}$.
For each element $K \in \CT_h$ and face $F \in \CF_h$, $\CV(K),\CV(F) \subset \CV_h$
are the sets of vertices of $K$ and $F$.

We require that the mesh is conforming in the sense of
\cite[Section 2.1.2]{boffi_brezzi_fortin_2013a} and
\cite[Definition 6.11]{ern_guermond_book}. Specifically, we assume that
$\overline{\Omega} = \cup_{K \in \CT_h} \overline{K}$, and that if
$K_\pm \in \CT_h$ are two distinct elements, $\overline{K_+} \cap \overline{K_-}$
is either empty, or a single vertex, edge or face of both $K_+$ and $K_-$.
This assumption is standard, and although it rules out hanging nodes, it does
not prevent heavily localized refinements with strong mesh grading \cite{apel_1999a}.

We further demand that the mesh fits the coefficient $\mm$, meaning that for all $K \in \CT_h$,
their exists $P \in \LP$ such that $\overline{K} \subset \overline{P}$. In effect, $\mm|_K$
is a constant value for all $K \in \CT_h$.

We finally assume that the surfaces $\{\Sigma_j\}_{j=1}^N$ characterizing the
space $\LBA_{\GN}(\Omega)$ are aligned with the mesh, meaning that for $1 \leq j \leq N$,
\begin{equation}
\label{eq_assumption_sigma}
\overline{\Sigma}_j \eq \bigcup_{F \in \CF_h^{\Sigma_j}} F
\end{equation}
for some $\CF_h^{\Sigma_j} \subset \CF_h$. Notice that because the particular choice
of surfaces $\{\Sigma_j\}_{j=1}^N$ is not important, but only their topological properties,
\eqref{eq_assumption_sigma} is by no means a restrictive assumption in practice
\cite{gross_kotiuga_2004a}. In fact, software packages are available to automatically
generate a set of surfaces $\{\Sigma_j\}_{j=1}^N$ aligned with faces, given
any conforming mesh as input, see e.g.
\cite{dlotko_kapidani_specogna_2017a,pellikka_suuriniemi_kettunen_geuzaine_2013a}.

\subsection{Mesh size and shape-regularity parameters}

For all $K \in \CT_h$, $h_K$ and $\rho_K$ respectively denote the diameters of
the smallest ball containing $K$ and the largest ball contained in $\overline{K}$.
We employ the notation $h_{\max} \eq \max_{K \in \CT_h} h_K$ for the mesh size.
The quantity $\kappa_K \eq h_K/\rho_K \geq 1$ is called the shape-regularity parameter
of $K$, and for $\CT \subset \CT_h$ we set $\kappa_{\CT} \eq \max_{K \in \CT} \kappa_K$.

\subsection{Jumps}

Each face $F \in \CF_h$ is equipped with a unit normal vector $\bn_F$.
If $F \in \CF_h^{\rm i}$ the orientation of $\bn_F$ is arbitrarily
fixed, and we assume that $\bn_F = \bn$ if $F \in \CF_h^{\rm e}$.
If $v: \Omega \to \mathbb R$ is a piecewise smooth function, its
jump through $F = \overline{K_+} \cap \overline{K_-} \in \CF_h^{\rm i}$
is defined by
\begin{equation*}
\jmp{v}_F
\eq
v_+|_F (\bn_+ \cdot \bn_F) + v_-|_F (\bn_- \cdot \bn_F)
\end{equation*}
where $v_\pm \eq v|_{K_\pm}$ and $\bn_\pm$ is the unit normal vector
of $K_\pm$. We then define the jump $\jmp{\bv}_F$ a piecewise smooth
vector-valued function $\bv: \Omega \to \mathbb R^3$ using the same
formula for each component.

\subsection{Vertex patches}
\label{section_vertex_patches}

For all $\ba \in \CV_h$, $\CTa \subset \CT_h$ collects the elements
$K \in \CT_h$ such that $\ba \in \CV(K)$. We also denote by $\oma$
the corresponding open domain. The hat function $\pa: \Omega \to \mathbb R$
is the only function affine in each element $K \in \CT_h$ such that
$\pa(\ba) = 1$ and $\pa(\bb) = 0$ for all $\bb \in \CV_h \setminus \{\ba\}$.
Notice that $\overline{\oma}$ then corresponds to the support of the hat function $\pa$.

We let $\gma \eq \GN \cup \{\pa = 0\}$ and
\begin{equation*}
\BH_0(\ddiv,\oma) \eq
\left \{
\bv \in \BH(\ddiv,\oma) \; | \;
\bv \cdot \bn = 0 \text{ on } \gma
\right \}.
\end{equation*}
We also introduce $L^2_\star(\oma) \eq \div \BH_0(\ddiv,\oma)$, so
that $L^2_\star(\oma)$ consists of functions with zero mean value
when $\gma = \partial \oma$, and coincide with the whole $L^2(\oma)$
otherwise. Finally, we set
\begin{equation*}
H^1_\star(\oma)
\eq
\left \{
v \in H^1(\oma) \; | \; v = 0 \text{ on } \partial \omega \setminus \gma
\right \}
\cap
L^2_\star(\oma).
\end{equation*}

The inequality
\begin{equation}
\label{eq_local_poincare_weighted}
\|\grad(\pa q)\|_{\mm,\oma}
\leq
\Cconta \LC_{\mm,\oma}^{1/2}\|\grad q\|_{\mm,\oma}
\qquad
\forall q \in H^1_\star(\oma)
\end{equation}
holds true with a constant $\Cconta$ only depending on the shape-regularity parameter
of $\CTa$, see, e.g., \cite[Equation (3.29)]{ern_vohralik_2015a}.

\subsection{Surfaces and half patches}

Consider a relatively open Lipschitz surface $\Sigma \subset \overline{\Omega}$
with $\partial \Sigma \subset \partial \Omega$ that corresponds
to a collection of faces of the mesh $\CT_h$. We assume that
$\Sigma$ is oriented by a unit normal vector $\bn_{\Sigma}$.
Consider a mesh vertex $\ba \in \CV_h$ such that $\ba \in \overline{\Sigma}$
and the associated patch $\oma$. Then we may uniquely decompose
$\oma$ into two connected subsets $\oma_{\Sigma,\pm}$ with outward unit
normal $\bn_{\oma_{\Sigma,\pm}}$ in such way that $\bn_{\oma_{\Sigma,\pm}} = \pm \bn_{\Sigma}$
on $\Sigma \cap \oma$. Notice that one of the two subsets may be empty if
$\Sigma \cap \oma \subset \partial \Omega$. These two halves of the patch
correspond to unions of mesh cells that we denote by $\CT_h^{\ba,\Sigma,\pm}$.
We also introduce the notation $\gma_{\Sigma,\pm} \eq \gma \cap \oma_{\Sigma,\pm}$.

Finally, we note that if $\ba \in \overline{\GN}$ and $\partial \Sigma \subset \overline{\GN}$,
then at least one of the two halves $\oma_{\Sigma,\pm}$ of the patch $\oma$ is not empty with the
boundary $\gma_{\Sigma,\pm}$ lying in $\GN$.

\subsection{Finite element spaces}

If $K \in \CT_h$ is an element of the mesh and $q \geq 0$,
then $\CP_q(K)$ is the set of polynomials mapping $K$ into $\mathbb R$
of degree at most $q$, and $\BCP_q(K) \eq [\CP_q(K)]^3$ collects vector-valued
polynomial functions. $\ND_q(K) \eq \BCP_q(K) \times \bx + \BCP_q(K)$
and $\RT_q(K) \eq \CP_q(K) \bx + \BCP_q(K)$ are then the space of N\'ed\'elec
\cite{nedelec_1980a} and Raviart--Thomas \cite{raviart_thomas_1977a} polynomials
of order $q$, see also \cite[Chapter 2.3]{boffi_brezzi_fortin_2013a} and
\cite[Chapters 11 and 12]{ern_guermond_book}.

For a subset of elements $\CT \subset \CT_h$, $\CP_q(\CT)$ stands for piecewise
polynomial functions on $\CT$, i.e. $v \in \CP_q(\CT)$ if and only if
$v|_K \in \CP_q(K)$ for all $K \in \CT$. $\BCP_q(\CT)$, $\ND_q(\CT)$
and $\RT_q(\CT)$ are defined similarly.

\subsection{Stable discrete minimization}

Following \cite[Theorem 2.3]{ern_vohralik_2021a} (see also
\cite[Proposition 3.1]{chaumontfrelet_vohralik_2023a} for general
configurations of boundary patches), for all vertices $\ba \in \CV_h$,
there exists a constant $\Csta$ only depending on the
shape-regularity parameter $\kappa_{\CTa}$ of $\CTa$ such that
\begin{equation}
\label{eq_stable_discrete_minimization}
\min_{\substack{
\bv_h \in \RT_q(\CTa) \cap \BH_0(\ddiv,\oma)
\\
\div \bv_h = r_h
}}
\|\bv_h-\bxi_h\|
\leq
\Csta
\min_{\substack{
\bv \in \BH_0(\ddiv,\oma)
\\
\div \bv = r_h
}}
\|\bv-\bxi_h\|
\end{equation}
for all polynomial degree $q \geq 0$, $r_h \in \CP_q(\CTa) \cap L^2_\star(\oma)$ and
$\bxi_h \in \RT_q(\CTa)$. Crucially, both minimizers in \eqref{eq_stable_discrete_minimization}
are uniquely defined, and $\Csta$ does not depend on $q$.

\subsection{Discrete solution}

Fix $p \geq 0$. We denote by $\BV_h \eq \ND_p(\CT_h) \cap \BH_{\GT}(\ccurl,\Omega)$
the usual space of N\'ed\'elec finite elements and set $\LBA_h \eq \curl \BV_h$.
There exists a unique couple $(\BH_h,\BA_h) \in \BV_h \times \LBA_h$ such that
\begin{equation}
\label{eq_magneto_discrete}
\left \{
\begin{array}{rcl}
(\mm\BH_h,\bv_h)_\Omega + (\BA_h,\curl\bv_h)_\Omega &=& 0,
\\
(\curl \BH_h,\bw_h)_\Omega &=& (\BJ,\bw_h)_\Omega.
\end{array}
\right .
\end{equation}
for all $\bv_h \in \BV_h$ and $\bw_h \in \LBA_h$.
Similar to \eqref{eq_magneto_continuous}, the existence and uniqueness of the solution
is due to the definition of $\LBA_h$ and the results in \cite{brezzi_1974a}.

\begin{remark}[Practical construction of $\LBA_h$]
Let us consider the set of Raviart--Thomas elements
$\BW_h \eq \RT_p(\CT_h) \cap \BH_{\GT}(\ddiv,\Omega)$.
Then, we have $\LBA_h \subset \BW_h$, and (at least)
two approaches are available to construct $\LBA_h$.

First (i), similar to \eqref{eq_condition_sigma}, it is known that
$\bw_h \in \LBA_h$ if and only if $\bw_h \in \BW_h$ with $\div \bw_h = 0$
and
\begin{equation}
\label{eq_condition_tilde_sigma}
\int_{\widetilde \Sigma_j} \bw_h \cdot \bn_{\widetilde \Sigma_j} = 0
\end{equation}
for $1 \leq j \leq \widetilde N$, where each $\widetilde \Sigma_j \subset \overline{\Omega}$
is a (relatively open) surface with $\partial \widetilde \Sigma_j \subset \GT$. We may assume
without generality that each $\widetilde \Sigma_j$ is exactly covered by mesh faces, so
that \eqref{eq_condition_tilde_sigma} may be efficiently checked numerically.
As previously mentioned, efficient algorithms are available to automatically construct
$\{\Sigma_j\}_{j=1}^{\widetilde N}$, see e.g.
\cite{%
dlotko_kapidani_specogna_2017a,%
gross_kotiuga_2004a,%
pellikka_suuriniemi_kettunen_geuzaine_2013a}.
Then, instead of explicitly constructing a basis of $\LBA_h$ we augment the discrete
formulation \eqref{eq_magneto_discrete} with two Lagrange multipliers, namely
$q_h \in \div \BW_h \subset \CP_p(\CT_h)$ to impose the divergence constraint,
and $\zeta \in \mathbb R^{\widetilde N}$ to enforce \eqref{eq_condition_tilde_sigma}.

As an alternative (ii), it is possible to employ a spanning-tree approach
to explicitly construct a basis of $\LBA_h$, see e.g.
\cite{alonsorodriguez_camano_delossantos_rapetti_2018a,alotto_perugia_1999a,scheichl_2002a}.
\end{remark}

\subsection{Oscillation term}

For $\bphi \in \BL^2(\Omega)$, consider the following orthogonal projector onto
the subset $\LBA_h$ of Raviart--Thomas elements
\begin{equation*}
\pi_h^0 \bphi
\eq
\arg \min_{\bv_h \in \LBA_h} \|\bphi-\bv_h\|_{\cc,\Omega}.
\end{equation*}
Then, the following quantity will be central for the data oscillation term
\begin{equation*}
\beta_h
\eq
\sup_{\substack{
\bphi \in \BH_{\GT}(\ccurl,\Omega) \cap \LBA_{\GN}(\Omega)
\\
\|\curl \bphi\|_{\cc,\Omega} = 1
}}
\|\bphi-\pi_h^0 \bphi\|_{\cc,\Omega}.
\end{equation*}

Notice that $\beta_h$ has dimension (length)$^{-1}$.
Actually, when the domain is convex (or smooth) and either $\GT$ or
$\GN$ is empty, the inclusion $\BH_{\GT}(\ccurl,\Omega) \cap \LBA_{\GN} \subset \BH^1(\Omega)$
holds true \cite[Theorem 2.17]{amrouche_bernardi_dauge_girault_1998a}, and $\beta_h \sim h_{\max}$.
A similar property still holds true in more general domains, if the mesh is properly graded
in the vicinity of re-entrant corners and edges \cite{apel_1999a,nicaise_2001a}.
In the general case, however, we have $\beta_h \sim \ell_\Omega^{1-s}h_{\max}^s$
for some $s > 0$ for which the inclusion
$\BH_{\GT}(\ccurl,\Omega) \cap \LBA_{\GN} \subset \BH^s(\Omega)$
holds true \cite{costabel_dauge_2000a}. When either $\GT$ or $\GN$ is empty, we can even select $s > 1/2$, see
\cite[Proposition 3.7]{amrouche_bernardi_dauge_girault_1998a}.

\section{Reliability}
\label{section_reliability}

\subsection{General upper bound}

We start with a general upper-bound given by any field $\BB \in \LBA_{\GN}(\Omega)$.
The bound \eqref{eq_prager_synge_magneto} is to be compared with the generalized
Prager--Synge bound \eqref{eq_prager_synge_electro} for the electrostatic problem.

\begin{theorem}[General upper-bound]
\label{theorem_upper_bound}
The estimate
\begin{equation}
\label{eq_prager_synge_magneto}
\|\BH-\BH_h\|_{\mm,\Omega}^2
\leq
\|\BB-\mm\BH_h\|_{\cc,\Omega}^2 + \beta_h^2 \|\BJ-\pi_h^0\BJ\|_{\mm,\Omega}^2
\end{equation}
holds true for all $\BB \in \LBA_{\GN}(\Omega)$.
\end{theorem}

\begin{proof}
Consider an arbitrary test function $\bv \in \BL^2(\Omega)$.
The space $\cc \LBA_{\GN}(\Omega)$ is closed in $\BL^2(\Omega)$,
so that we may introduce the decomposition $\bv = \bv_0 + \bv^\perp$,
with $\bv_0 \in \cc \LBA_{\GN}(\Omega)$ and
$\bv^\perp \in \left (\cc\LBA_{\GN}(\Omega)\right )^\perp$,
orthogonal for the $(\mm\cdot,\cdot)_\Omega$ inner-product.
Then, we have
\begin{equation*}
(\mm(\BH-\BH_h),\bv)_\Omega
=
(\BH-\BH_h,\mm\bv_0)_\Omega + (\mm(\BH-\BH_h),\bv^\perp)_\Omega.
\end{equation*}

The first equation of the magnetostatic problem \eqref{eq_magneto_continuous}
ensures that $\mm \BH \in \LBA_{\GN}(\Omega)$. Hence, $\BH \in \cc \LBA_{\GN}(\Omega)$,
and $(\mm\BH,\bv^\perp)_\Omega = 0$. As a result, we have
\begin{equation*}
(\mm(\BH-\BH_h),\bv^\perp)_\Omega = (\mm(\cc\BB-\BH_h),\bv^\perp)_\Omega
\end{equation*}
for all $\BB \in \LBA_{\GN}(\Omega)$.

On the other hand, since $\bv_0 \in \cc\LBA_{\GN}(\Omega)$, there exists
$\bphi \in \BH_{\GN}(\ccurl,\Omega)$ such that $\mm \bv_0 = \curl \bphi$.
In addition, since only $\curl \bphi$ (as opposed to $\bphi$ itself) intervenes,
we may select $\bphi$ such that $\bphi \in \LBA_{\GT}(\Omega)$. Specifically,
this may be done by selecting $\bpsi \in \BH_{\GT}(\ccurl,\Omega)$ such that
$\curl \curl \bpsi = \bv_0$ in $\Omega$, e.g. by solving \eqref{eq_magneto_continuous}
with $\bv_0$ instead of $\BJ$ as a right-hand side, and then letting $\bphi = \curl \bpsi$.
Then, we have
\begin{equation*}
(\BH-\BH_h,\mm\bv_0)_\Omega
=
(\BH-\BH_h,\curl \bphi)_\Omega
=
(\BJ-\pi_h^0 \BJ,\bphi)_\Omega.
\end{equation*}
Since $\pi_h^0$ is an orthogonal projector, we can then write that
\begin{equation*}
(\BH-\BH_h,\mm\bv_0)_\Omega
=
(\BJ-\pi_h^0 \BJ,\bphi-\pi_h^0 \bphi)_\Omega
\leq
\|\BJ-\pi_h^0\BJ\|_{\mm,\Omega}\|\bphi-\pi_h^0\bphi\|_{\cc,\Omega},
\end{equation*}
and it follows from the definition of $\beta_h$ and the fact
that $\curl \bphi = \mm\bv_0$ that
\begin{equation*}
(\BH-\BH_h,\mm\bv_0)_\Omega
\leq
\beta_h \|\BJ-\pi_h^0\BJ\|_{\mm,\Omega}\|\curl \bphi\|_{\cc,\Omega}
=
\beta_h \|\BJ-\pi_h^0\BJ\|_{\mm,\Omega}\|\bv_0\|_{\mm,\Omega}.
\end{equation*}
We can then conclude the proof with
\begin{align*}
(\mm(\BH-\BH_h),\bv)
&\leq
\beta_h\|\BJ-\pi_h^0\BJ\|_{\mm,\Omega}\|\bv_0\|_{\mm,\Omega}
+
\|\cc\BB-\BH_h\|_{\mm,\Omega}\|\bv^\perp\|_{\mm,\Omega},
\\
&\leq
\left (
\beta_h^2\|\BJ-\pi_h^0\BJ\|_{\mm,\Omega}^2
+
\|\cc\BB-\BH_h\|_{\mm,\Omega}^2
\right )^{1/2}
\left (
\|\bv_0\|_{\mm,\Omega}^2
+
\|\bv^\perp\|_{\mm,\Omega}^2
\right )^{1/2}
\\
&=
\left (
\beta_h^2\|\BJ-\pi_h^0\BJ\|_{\mm,\Omega}^2
+
\|\cc\BB-\BH_h\|_{\mm,\Omega}^2
\right )^{1/2}
\|\bv\|_{\mm,\Omega},
\end{align*}
since $\bv$ was arbitrary in $\BL^2(\Omega)$.
\end{proof}

\subsection{Magnetic field reconstruction}
\label{section_reconstruction}

We now present a practical construction of a discrete field $\BB_h \in \LBA_{\GN}(\Omega)$
that is obtained as local post-processing of $\BH_h$. Our construction is based on
divergence-constrained patch-wise minimization problems, and in fact closely
follows the equilibrated-flux reconstruction of $\BD_h$ in electrostatic problems
\cite{braess_pillwein_schoberl_2009a,destuynder_metivet_1999a,ern_vohralik_2015a}.

For each vertex $\ba \in \CV_h$, the local contribution to our induction field
reconstruction is given by
\begin{equation}
\label{eq_Bh_local}
\BB_h^{\ba}
\eq
\arg \min_{\substack{
\bv_h \in \RT_{p+2}(\CTa) \cap \BH_0(\ddiv,\oma)
\\
\div \bv_h = \grad \pa \cdot (\mm\BH_h)
}}
\|\bv_h-\pa(\mm\BH_h)\|_{\cc,\oma}.
\end{equation}
This definition indeed makes sense: the compatibility condition is satisfied
due to the first equation in \eqref{eq_magneto_discrete}, since
$\grad \pa \in \ND_0(\CT_h) \cap \BH_0(\ccurl,\Omega)$. Introducing
\begin{equation}
\label{eq_Bh_global}
\BB_h \eq \sum_{\ba \in \CV_h} \BB_h^{\ba},
\end{equation}
it is then clear that $\BB_h \in \BH_{\GN}(\ddiv,\Omega)$ with $\div \BB_h = 0$.

\begin{remark}[Polynomial degree of the reconstruction]
The polynomial degree of the reconstruction is increased by two as compared
to the finite element solution. This may seem surprising, as in the case
of electrostatic problems, it is sufficient to increase the polynomial degree
by one \cite{braess_pillwein_schoberl_2009a,destuynder_metivet_1999a,ern_vohralik_2015a}.
The reason for increasing to polynomial degree by one in electrostatic
problems is because of the multiplication by the hat function $\pa$. Here, we
have the additional issue that we are trying to represent a (broken) N\'ed\'elec
polynomial with Raviart--Thomas polynomials, which comes at the expense of
another polynomial degree increase. It is possible to only increase the polynomial
degree by one if Brezzi--Douglas--Marini are used instead of Raviart--Thomas elements
in the reconstruction, but we do not pursue this approach here for the sake of simplicity.
\end{remark}

For configurations where cohomology spaces are trivial, the above considerations
suffice to infer that $\BB_h \in \LBA_h$. In the general case however, we must
additionally check that the periods of $\BB_h$ vanish according to \eqref{eq_condition_sigma}.
Surprisingly, this is automatically satisfied, as we next establish.

\begin{lemma}[Periods of the reconstruction]
\label{lemma_condition_sigma}
Let $\Sigma \subset \Omega$ be a relatively open Lipschitz surface consisting
of a collection of mesh faces such that $\partial \Sigma \subset \GN$. If $\Sigma$
if oriented with a unit normal vector $\bn_{\Sigma}$, then, we have
\begin{equation*}
\int_{\Sigma} \BB_h \cdot \bn_{\Sigma} = 0.
\end{equation*}
\end{lemma}

\begin{proof}
Let us denote by $\CV_h^{\Sigma}$ the set of vertices $\ba \in \CV_h$
belonging to $\overline{\Sigma}$. For each $\ba \in \CV_h^{\Sigma}$,
we let $\iota_{\ba} \in \{-,+\}$ be such that $\oma_{\Sigma,\iota_{\ba}}$
is non-empty and, if $\ba \in \overline{\GN}$, such that
$\gma_{\Sigma,\iota_{\ba}} \subset \overline{\GN}$.

It is then clear that for each $\ba \in \CV_h^{\Sigma}$
\begin{multline*}
\int_{\Sigma} \BB_h^{\ba} \cdot \bn_{\Sigma}
=
\iota_{\ba} \int_{\partial \oma_{\Sigma,\iota_{\ba}}}
\BB_h^{\ba} \cdot \bn_{\oma_{\Sigma,\iota{\ba}}}
=
\iota_{\ba} \int_{\oma_{\Sigma,\iota_{\ba}}} \div \BB_h^{\ba}
\\
=
\int_{\oma_{\Sigma,\iota_{\ba}}}
(\iota_{\ba} \grad \pa) \cdot (\mm\BH_h)
=
(\mm\BH_h, \iota_{\ba} \chi_{\ba} \grad \pa)_\Omega
\end{multline*}
where we let $\chi_{\ba} \in \CP_0(\CT_h)$ denote the set function
of $\oma_{\Sigma,\iota_{\ba}}$. As a result, introducing the test function
\begin{equation}
\label{tmp_def_test_function}
\bv_h \eq \sum_{\ba \in \CV_h^{\Sigma}} \iota_{\ba} \chi_{\ba} \grad \pa \in \BL^2(\Omega),
\end{equation}
we have
\begin{equation*}
\int_{\Sigma} \BB_h^{\ba} \cdot \bn_{\Sigma} = (\mm\BH_h,\bv_h)_\Omega,
\end{equation*}
and using the first equation of \eqref{eq_magneto_discrete},
the conclusion follows if we can show that $\bv_h \in \BV_h$
with $\curl \bv_h = \bo$.

We first observe that in each cell $K \in \CT_h$, $(\iota_{\ba}\chi_{\ba})|_K$
is constant, so that
\begin{equation*}
\bv_h|_K
=
\sum_{\ba \in \CV_h^{\Sigma}} (\iota_{\ba}\chi_{\ba})|_K \grad (\pa|_K)
=
\grad \widetilde q_K
\end{equation*}
for some $\widetilde q_K \in \CP_1(K)$. As a result
$\bv_h \in \N_0(\CT_h)$ with $\curl \bv_h = \bo$ element-wise,
and therefore, the conclusion follows if we can show that
$\jmp{\bv_h}_F \times \bn_F = \bo$ for all $F \in \CF_h^{\rm i}$
and $\bv_h \times \bn = \bo$ on $\GT$.

Let us first consider an interior face $F \in \CF_h^{\rm i}$.
Notice that then $\jmp{\grad \pa}_F \times \bn_F = \bo$
for all $\ba \in \CV_h$, and we have
\begin{equation*}
\jmp{\bv_h}_F \times \bn_F
=
\sum_{\ba \in \CV_h^{\Sigma}}
\iota_{\ba} \jmp{\chi_{\ba}}_F \grad \pa \times \bn_F
=
\sum_{\ba \in \CV(F)}
\iota_{\ba} \jmp{\chi_{\ba}}_F \grad \pa \times \bn_F.
\end{equation*}
We then distinguish two cases. (i) If $F \not \subset \overline{\Sigma}$,
for each $\ba \in \CV(F)$, then either (ia) $F \subset \partial \oma$,
in which case, $\grad \pa \times \bn_F = \bo$ and $\jmp{\bv_h}_F \times \bn_F = \bo$
or (ib) $F \not\subset\partial \oma$, in which case $\jmp{\chi_{\ba}} = 0$,
again leading to $\jmp{\bv_h}_F \times \bn_F = \bo$. On the other hand, (ii) if
$F \subset \Sigma$, we have
\begin{equation*}
\jmp{\bv_h}_F \times \bn_F
=
\sum_{\ba \in \CV(F)} \grad \pa \times \bn_F
\end{equation*}
since $\jmp{\chi_{\ba}}_F = \iota_{\ba}$ for all the relevant vertices in the sum.
But then, we easily conclude since
\begin{equation*}
\jmp{\bv_h}_F \times \bn_F
=
\grad \left (\sum_{\ba \in \CV(F)} \pa\right ) \times \bn_F
=
\grad \left (1\right ) \times \bn_F
=
\bzero.
\end{equation*}

We conlude the proof by considering a face $F \subset \overline{\GT}$,
whereby we need to show that $\bv_h \times \bn_F = \bzero$ on $F$. We first observe
that due to the locality of the hat functions $\pa$, the value
of $\bv_h \times \bn_F$ on $F$ can only be influenced by
boundary vertices $\ba \subset \partial \Omega$ in \eqref{tmp_def_test_function}.
Since on the other hand $\partial \Sigma \subset \overline{\GN}$,
such vertices must lie in $\overline{\GN}$. However, for such vertices,
we have always chosen the half-patch $\oma_{\Sigma,\iota_{\ba}}$
such that $\gma_{\Sigma,\iota_{\ba}} \subset \overline{\GN}$, meaning
that corresponding contribution vanishes on $\GT$, as $\chi_{\ba} = 0$ on $\GT$.
\end{proof}

A direct consequence of the characterization \eqref{eq_condition_sigma} of the image
of the curl, the assumption \eqref{eq_assumption_sigma} that the surfaces $\{\Sigma_j\}_{j=1}^N$
are aligned with the mesh $\CT_h$ and Lemma \ref{lemma_condition_sigma} above is that $\BB_h$ is
indeed the curl of some potential.

\begin{theorem}[Localized reconstruction]
For the induction field reconstruction provided by \eqref{eq_Bh_local}
and \eqref{eq_Bh_global}, we have $\BB_h \in \LBA_{\GN}(\Omega)$.
\end{theorem}

\section{Efficiency}
\label{section_efficiency}

In this section, we present efficiency results for the estimator constructed
with the procedure of Section \ref{section_reconstruction}. We start with a
patch-wise efficiency result that applies to each local contribution $\BB_h^{\ba}$.
The proof is similar to seminal results for the electrostatic problem, see
\cite{braess_schoberl_2008a,ern_vohralik_2015a}.

\begin{theorem}[Lower-bound]
The estimate
\begin{equation}
\label{eq_lower_bound_patch}
\|\BB_h^{\ba}-\pa(\mm\BH_h)\|_{\cc,\oma}
\leq
\Csta \Cconta \LC_{\mm,\oma}
\|\BH-\BH_h\|_{\mm,\oma}
\end{equation}
holds true for all vertices $\ba \in \CV_h$.
\end{theorem}

\begin{proof}
Let $\BB^{\ba}$ solve the continuous version of \eqref{eq_Bh_local}, i.e.
\begin{equation}
\label{tmp_B_local}
\BB^{\ba}
\eq
\arg \min_{\substack{
\bv \in \BH_0(\ddiv,\oma)
\\
\div \bv = \grad \pa \cdot (\mm\BH_h)
}}
\|\bv-\pa(\mm\BH_h)\|_{\cc,\oma}.
\end{equation}
Then, \eqref{eq_stable_discrete_minimization} ensures that
\begin{equation}
\|\BB_h^{\ba}-\pa(\mm\BH_h)\|_{\cc,\oma}
\leq
\Csta
\LC_{\mm,\oma}^{1/2}
\|\BB^{\ba}-\pa(\mm\BH_h)\|_{\cc,\oma}.
\end{equation}
Indeed, $\BH_h \in \ND_p(\CT_h) \subset \BCP_{p+1}(\CT_h)$, so that
$\pa (\mm\BH_h) \in \BCP_{p+2}(\CT_h) \subset \RT_{p+2}(\CT_h)$.
As a result, we only need to focus on the continuous minimization problem \eqref{tmp_B_local}.
Notice that the additional factor $\LC_{\mm,\oma}^{1/2}$ as compared to
\eqref{eq_stable_discrete_minimization} is due to the weighted norms.

The Euler--Lagrange equations for \eqref{tmp_B_local} show that there exists
a unique $r^{\ba} \in L^2_\star(\oma)$ such that
\begin{equation*}
\left \{
\begin{aligned}
(\cc \BB^{\ba},\bv)_{\oma} + (r^{\ba},\div \bv)_{\oma}
&=
(\pa \BH_h,\bv)_{\oma}
\\
(\div \BB^{\ba},q)_{\oma}
&=
(\grad \pa \cdot (\mm\BH_h),q)_{\oma}
\end{aligned}
\right .
\end{equation*}
for all $\bv \in \BH_0(\ddiv,\oma)$ and $q \in L^2_\star(\oma)$.
Integrating by parts in the first equation shows that we actually have
$r^{\ba} \in H^1_\star(\oma)$ with
\begin{equation*}
\mm \grad r^{\ba} = \BB^{\ba}-\pa(\mm\BH_h),
\end{equation*}
and in particular
\begin{equation*}
\|\grad r^{\ba}\|_{\mm,\oma}
=
\|\cc(\BB^{\ba}-\pa(\mm\BH_h))\|_{\mm,\oma}
=
\|\BB^{\ba}-\pa(\mm\BH_h)\|_{\cc,\oma}.
\end{equation*}
To complete the proof, we estimate the norm of $\grad r^{\ba}$. We start with
\begin{align*}
(\mm \grad r^{\ba},\grad v)_{\oma}
&=
(\BB^{\ba},\grad v)_{\oma} - (\pa(\mm\BH_h),\grad v)_{\oma}
\\
&=
-(\div \BB^{\ba},v)_{\oma} - (\pa(\mm\BH_h),\grad v)_{\oma}
\\
&=
-(\grad \pa \cdot (\mm\BH_h),v)_{\oma} - (\pa(\mm\BH_h),\grad v)_{\oma}
\\
&=
-(\mm\BH_h,\grad(\pa v))_{\oma},
\end{align*}
and upon recalling that $\div(\mm\BH) = 0$, we arrive at
\begin{equation*}
(\mm \grad r^{\ba},\grad v)_{\oma}
=
(\mm(\BH-\BH_h),\grad(\pa v))_{\oma}.
\end{equation*}
At that point, we easily finish the proof with \eqref{eq_local_poincare_weighted} since
\begin{align*}
\|\grad r^{\ba}\|_{\mm,\oma}^2
\leq
\|\BH-\BH_h\|_{\mm,\oma}\|\grad(\pa r^{\ba})\|_{\mm,\oma}
\leq
\Cconta
\LC_{\mm,\oma}^{1/2}
\|\BH-\BH_h\|_{\mm,\oma}
\|\grad r^{\ba}\|_{\mm,\oma},
\end{align*}
so that
\begin{equation*}
\|\grad r^{\ba}\|_{\mm,\oma}
\leq
\Cconta
\LC_{\mm,\oma}^{1/2}
\|\BH-\BH_h\|_{\mm,\oma}.
\end{equation*}
\end{proof}

From the efficiency of the local contribution $\BB_h^{\ba}$, we finally
establish an efficiency result for the induction field reconstruction $\BB_h$.

\begin{corollary}[Lower-bounds for the element-wise estimator]
\label{corollary_lower_bound}
The efficiency estimate
\begin{equation}
\label{eq_lower_bound_element}
\|\BB_h-\mm\BH_h\|_{\cc,K}
\leq
2\CstK \CcontK \LC_{\mm,\tK}
\|\BH-\BH_h\|_{\mm,\tK}
\end{equation}
holds true for all $K \in \CT_h$ with
\begin{equation*}
\tK \eq \bigcup_{\ba \in \CV(K)} \oma,
\qquad
\CstK \eq \max_{\ba \in \CV(K)} \Csta,
\qquad
\CcontK \eq \max_{\ba \in \CV(K)} \CcontK.
\end{equation*}
\end{corollary}

\begin{proof}
Let $K \in \CT_h$. The four hat functions associated with the vertices of $K$
form a partition of unity over $K$, and $\BB_h^{\ba} = \bo$ on $K$ unless
$\ba \in \CV(K)$. As a result, we have
\begin{equation*}
\BB_h-\mm\BH_h
=
\sum_{\ba \in \CV(K)}
\BB^{\ba}_h-\pa(\mm\BH_h)
\end{equation*}
on $K$, and the triangle inequality immediately gives that
\begin{align*}
\|\BB_h-\mm\BH_h\|_{\cc,K}
=
\left \|
\sum_{\ba \in \CV(K)}
\BB^{\ba}_h-\pa(\mm\BH_h)
\right \|_{\cc,K}
&\leq
\sum_{\ba \in \CV(K)}
\left \|
\BB^{\ba}_h-\pa(\mm\BH_h)
\right \|_{\cc,K}
\\
&\leq
\sum_{\ba \in \CV(K)}
\left \|
\BB^{\ba}_h-\pa(\mm\BH_h)
\right \|_{\cc,\oma}.
\end{align*}
The estimate in \eqref{eq_lower_bound_element} then follows from
\eqref{eq_lower_bound_patch} since
\begin{align*}
\|\BB_h-\mm\BH_h\|_{\cc,K}
&\leq
\sum_{\ba \in \CV(K)}
\Csta \Cconta \LC_{\mm,\oma}
\|\BH_h-\BH_h\|_{\mm,\oma}
\\
&\leq
\CstK \CcontK \LC_{\mm,\tK}
\sum_{\ba \in \CV(K)}
\|\BH_h-\BH_h\|_{\mm,\oma}
\\
&\leq
2 \CstK \CcontK \LC_{\mm,\tK}
\|\BH_h-\BH_h\|_{\mm,\tK},
\end{align*}
as $K$ has $4$ vertices.
\end{proof}

\section{Numerical examples}
\label{section_numerical_examples}

\subsection{Settings}

In all numerical examples, we consider the case where $\GN = \emptyset$, and thus,
$\GT = \partial \Omega$. For simplicity, we also set $\mm \equiv \BI$. Notice
that then $\BH = \BB$.

In the graphs below, we employ the notations
\begin{equation*}
\merrH \eq \|\BH-\BH_h\|_\Omega,
\qquad
\merrB \eq \|\BH-\BB_h\|_\Omega,
\qquad
\mest \eq \|\BH_h-\BB_h\|_\Omega,
\end{equation*}
where $\BH_h$ is the finite element solution, and $\BB_h$ is
the reconstruction introduced in Section \ref{section_reconstruction}.
Observe that for the sake of simplicity, we have not introduced
the oscillation term in $\mest$. Element-wise versions of $\merrH$
and $\mest$ are also plotted in the adaptivity example. The notation
$N_{\rm dofs}$ stands for the dimension of the N\'ed\'elec finite element
space $\BV_h$.

We employ {\tt gmsh} to generate meshes \cite{geuzaine_remacle_2009a}, and {\tt mmg3d}
for adaptive mesh refinements \cite{dobrzynski_2012a}. The linear systems are solved
using the {\tt mumps} software package \cite{amestoy_duff_lexcellent_2000a}.

\subsection{Trivial cohomology}

Here, $\Omega = (0,1)^3$ is the unit cube. The right-hand side and solution respectively read
\begin{equation}
\label{eq_rhs_solution}
\BJ
\eq
\left (
\begin{array}{c}
\sin(\pi\bx_1)\cos(\pi\bx_2)
\\
-\cos(\pi\bx_1)\sin(\pi\bx_2)
\\
0
\end{array}
\right ),
\qquad
\BH
\eq
\left (
\begin{array}{c}
0
\\
0
\\
\sin(\pi\bx_1)\sin(\pi\bx_2)
\end{array}
\right ).
\end{equation}

\begin{figure}
\begin{minipage}{.45\linewidth}
\begin{tikzpicture}
\begin{axis}
[
width=\textwidth,
ymode = log,
xmode = log,
xlabel={\ndof}
]

\plot[thick,blue ,mark=o     ] table[x=dofs,y=errH] {data/trivial/P1.txt}
node[pos=0.6,pin={[pin distance=0.1cm]-90:{\errH}}] {};
\plot[thick,red  ,mark=x     ] table[x=dofs,y=errB] {data/trivial/P1.txt}
node[pos=0.1,pin={[pin distance=0.8cm]-90:{\errB}}] {};
\plot[thick,black,mark=square] table[x=dofs,y=est ] {data/trivial/P1.txt}
node[pos=0.9,pin={[pin distance=0.5cm] 90:{\est }}] {};

\plot [dashed,domain=1e3:5e5] {0.3*x^(-2./3.)};

\SlopeTriangle{0.55}{-0.2}{0.15}{-0.66666}{$\ndofs^{-2/3}$}{}

\end{axis}
\end{tikzpicture}
\end{minipage}
\begin{minipage}{.45\linewidth}
\begin{tikzpicture}
\begin{axis}
[
width=\textwidth,
xmode = log,
ymin = 0.9,
ymax = 1.1,
xlabel={\ndof}
]

\plot[thick,black,mark=square] table[x=dofs,y expr=\thisrow{est}/\thisrow{errH}]
{data/trivial/P1.txt}
node[pos=0.5,pin={[pin distance=0.5cm]-90:{\est/\errH}}] {};

\end{axis}
\end{tikzpicture}
\end{minipage}
\caption{Trivial cohomology example, $p=1$}
\label{figure_trivial_P1}
\end{figure}
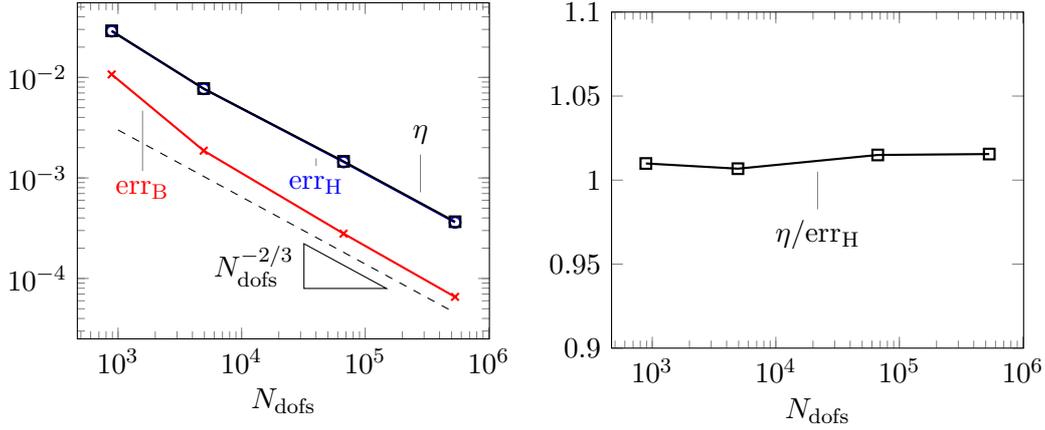

In Figures \ref{figure_trivial_P0} and \ref{figure_trivial_P1},
we respectively fix the polynomial degree to $p=0$ and $1$,
and let the mesh size vary. The expected convergence rates (namely
linear and quadratic) are observed for the discrete solution
$\BH_h$, and the estimator closely follows the error. Indeed,
as can be seen on the right panels of Figures \ref{figure_trivial_P0}
and \ref{figure_trivial_P1}, the effectivity indices stay slightly
above one. We also see that the reconstruction $\BB_h$ converges
toward $\BB$ with the same rates as $\BH_h$ converges to $\BH$.
The post-processing is about $3$ to $5$ times more accurate than
the original discrete solution in this case.

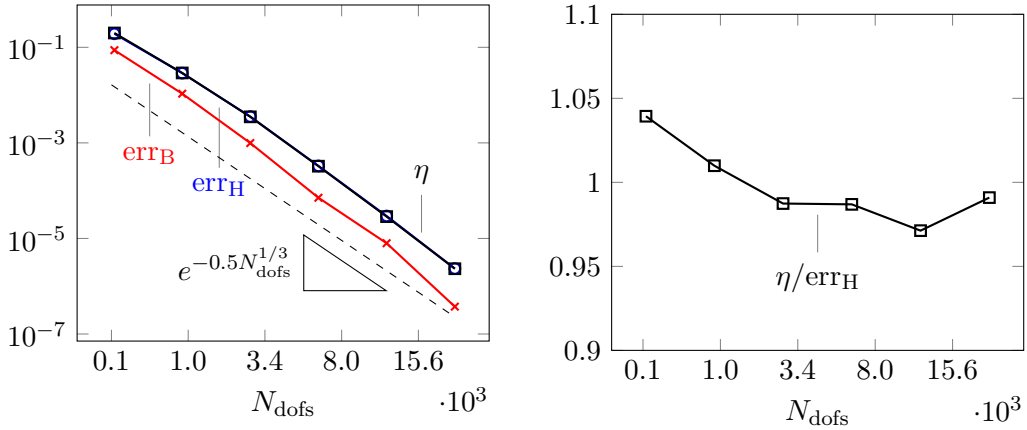
\begin{figure}
\begin{minipage}{.45\linewidth}
\begin{tikzpicture}
\begin{axis}
[
width=\textwidth,
ymode = log,
xlabel={\ndof},
xtick={5,10,15,20,25},
xticklabels={0.1,1.0,3.4,8.0,15.6},
scaled x ticks={real:1.0},
xtick scale label code/.code={$\cdot 10^3$}]
]

\plot[thick,blue ,mark=o     ] table[x expr=\thisrow{dofs}^(1./3.),y=errH]
{data/trivial/H.txt} node[pos=0.3,pin={[pin distance=0.8cm]-90:{\errH}}] {};
\plot[thick,red  ,mark=x     ] table[x expr=\thisrow{dofs}^(1./3.),y=errB]
{data/trivial/H.txt} node[pos=0.1,pin={[pin distance=0.7cm]-90:{\errB}}] {};
\plot[thick,black,mark=square] table[x expr=\thisrow{dofs}^(1./3.),y=est ]
{data/trivial/H.txt} node[pos=0.9,pin={[pin distance=0.5cm] 90:{\est }}] {};

\plot [dashed,domain=5:27] {0.2*exp(-0.5*x)};

\SlopeTriangle{0.55}{-0.2}{0.15}{-0.5}{$e^{-0.5\ndofs^{1/3}}$}{}
\end{axis}
\end{tikzpicture}
\end{minipage}
\begin{minipage}{.45\linewidth}
\begin{tikzpicture}
\begin{axis}
[
width=\textwidth,
ymin = 0.9,
ymax = 1.1,
xlabel={\ndof},
xtick={5,10,15,20,25},
xticklabels={0.1,1.0,3.4,8.0,15.6},
scaled x ticks={real:1.0},
xtick scale label code/.code={$\cdot 10^3$}]
]

\plot[thick,black,mark=square]
table[x expr=\thisrow{dofs}^(1./3.),y expr=\thisrow{est}/\thisrow{errH}]
{data/trivial/H.txt}
node[pos=0.5,pin={[pin distance=0.5cm]-90:{\est/\errH}}] {};

\end{axis}
\end{tikzpicture}
\end{minipage}
\caption{Trivial cohomology example, $p$-convergence}
\label{figure_trivial_H}
\end{figure}

Figure \ref{figure_trivial_H} displays a $p$-convergence result
whereby a mesh of $176$ tetrahedra is fixed ($h_{\max} = 0.76$),
and $p$ varies from $0$ to $5$. On the left panel, we note that the
convergence is exponential in the number of dofs, as is to be expected.
The post-processed solution $\BB_h$ is $2$ to $4$ times more accurate than
$\BH_h$. The behaviour of the effectivity index is reported on the
right-panel. It slightly drops below one because the mesh is fairly
coarse and the data oscillation term have not been computed. Yet,
the agreement between the error and the estimator is excellent, and
the effectivity index remains essentially constant for all polynomial
degrees, highlighting the $p$-robustness of the estimator.

\begin{figure}
\begin{minipage}{.45\linewidth}
\begin{tikzpicture}
\begin{axis}
[
width=\textwidth,
ymode = log,
xmode = log,
xlabel={\ndof}
]

\plot[thick,blue ,mark=o     ] table[x=dofs,y=errH] {data/trivial/P0.txt}
node[pos=0.5,pin={[pin distance=0.15cm]-90:{\errH}}] {};
\plot[thick,red  ,mark=x     ] table[x=dofs,y=errB] {data/trivial/P0.txt}
node[pos=0.1,pin={[pin distance=0.7cm]-90:{\errB}}] {};
\plot[thick,black,mark=square] table[x=dofs,y=est ] {data/trivial/P0.txt}
node[pos=0.9,pin={[pin distance=0.5cm] 90:{\est }}] {};

\plot [dashed,domain=1e2:1e5] {0.25*x^(-1./3.)};

\SlopeTriangle{0.55}{-0.2}{0.15}{-0.333333}{$\ndofs^{-1/3}$}{}

\end{axis}
\end{tikzpicture}
\end{minipage}
\begin{minipage}{.45\linewidth}
\begin{tikzpicture}
\begin{axis}
[
width=\textwidth,
xmode = log,
ymin = 0.9,
ymax = 1.1,
xlabel={\ndof}
]

\plot[thick,black,mark=square] table[x=dofs,y expr=\thisrow{est}/\thisrow{errH}]
{data/trivial/P0.txt}
node[pos=0.5,pin={[pin distance=0.5cm]-90:{\est/\errH}}] {};

\end{axis}
\end{tikzpicture}
\end{minipage}
\caption{Trivial cohomology, $p=0$}
\label{figure_trivial_P0}
\end{figure}
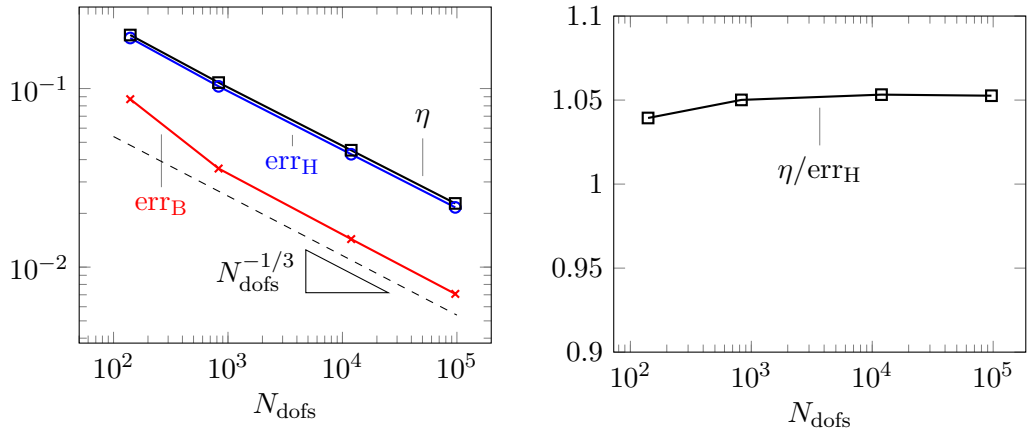

\subsection{Non-trivial cohomology}

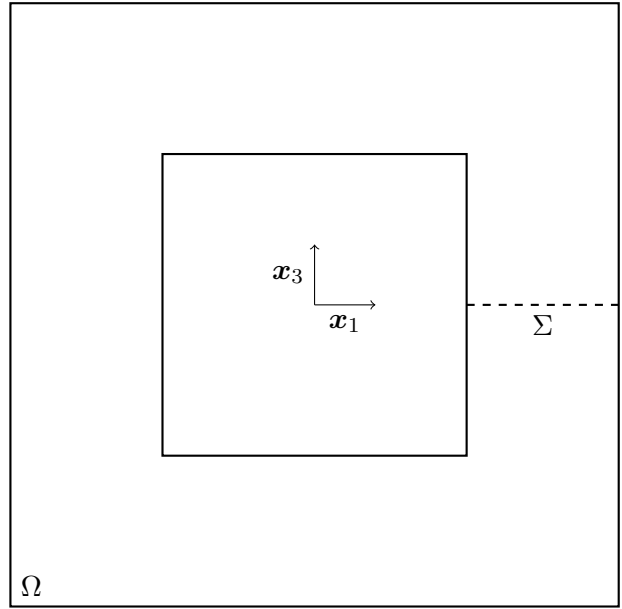
\begin{figure}
\begin{tikzpicture}[scale=2]
\draw[thick] (-2,-2) rectangle (2,2);
\draw[thick] (-1,-1) rectangle (1,1);
\draw[thick,dashed] ( 1, 0) -- (2,0);
\draw (1.5,0) node[anchor=north] {$\Sigma$};

\draw[->] (0.0,0.0) -- (0.4,0.0);
\draw[->] (0.0,0.0) -- (0.0,0.4);

\draw (0.2,0.0) node[anchor=north] {$\bx_1$};
\draw (0.0,0.2) node[anchor=east]  {$\bx_3$};

\draw(-2,-2) node[anchor=south west] {$\Omega$};
\end{tikzpicture}
\caption{Geometry of the non-trivial topology example: cut at $\bx_2 = 0$}
\label{figure_nontrivial_sketch}
\end{figure}

We then consider an example where the cohomology space is non-trivial  by setting
$\Omega \eq R \setminus C$ with $R \eq (-2,2) \times (-1,1) \times (-2,2)$ and
$C \eq (-1,1)^3$. $\Omega$ thus has the topology of a torus. To easily
construct $\Lambda_h$, we require that the surface $\Sigma \eq (1,2) \times (-1,1) \times \{0\}$
is exactly meshed, and impose that $(\bw_h \cdot \bn,1)_{\Sigma} = 0$ for all
$\bw_h \in \Lambda_h$. Figures \ref{figure_nontrivial_sketch} sketches the
geometry of the domain. The right-hand side and solution are still defined
with \eqref{eq_rhs_solution}.

\begin{figure}
\begin{minipage}{.45\linewidth}
\begin{tikzpicture}
\begin{axis}
[
width=\textwidth,
ymode = log,
xmode = log,
xlabel={\ndof}
]

\plot[thick,blue ,mark=o     ] table[x=dofs,y=errH] {data/nontrivial/P0.txt}
node[pos=0.5,pin={[pin distance=0.2cm]-90:{\errH}}] {};
\plot[thick,red  ,mark=x     ] table[x=dofs,y=errB] {data/nontrivial/P0.txt}
node[pos=0.1,pin={[pin distance=0.7cm]-90:{\errB}}] {};
\plot[thick,black,mark=square] table[x=dofs,y=est ] {data/nontrivial/P0.txt}
node[pos=0.9,pin={[pin distance=0.5cm] 90:{\est }}] {};

\plot [dashed,domain=1e3:1e5] {5*x^(-1./3.)};

\SlopeTriangle{0.55}{-0.2}{0.1}{-0.333333}{$\ndofs^{-1/3}$}{}

\end{axis}
\end{tikzpicture}
\end{minipage}
\begin{minipage}{.45\linewidth}
\begin{tikzpicture}
\begin{axis}
[
width=\textwidth,
xmode = log,
ymin = 0.9,
ymax = 1.1,
xlabel={\ndof}
]

\plot[thick,black,mark=square] table[x=dofs,y expr=\thisrow{est}/\thisrow{errH}]
{data/nontrivial/P0.txt}
node[pos=0.5,pin={[pin distance=0.5cm]-90:{\est/\errH}}] {};

\end{axis}
\end{tikzpicture}
\end{minipage}
\caption{Non-trivial cohomology example, $p=0$}
\label{figure_nontrivial_P0}
\end{figure}
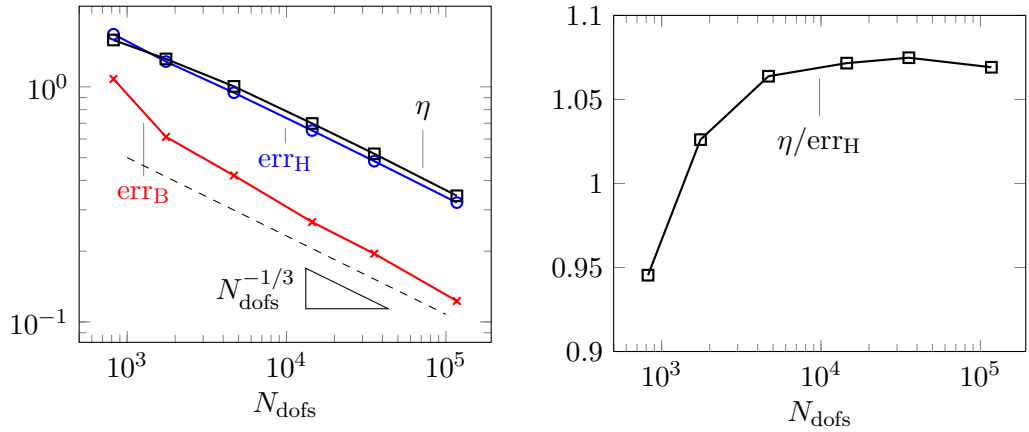

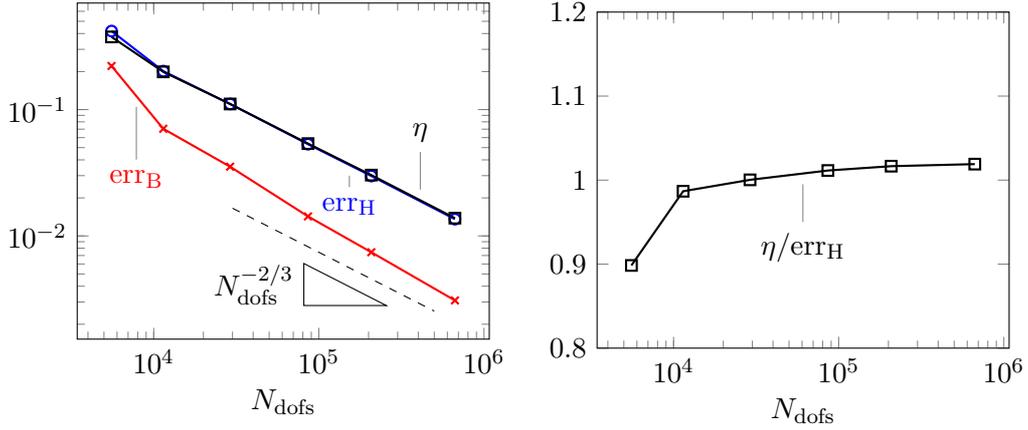
\begin{figure}
\begin{minipage}{.45\linewidth}
\begin{tikzpicture}
\begin{axis}
[
width=\textwidth,
ymode = log,
xmode = log,
xlabel={\ndof}
]

\plot[thick,blue ,mark=o     ] table[x=dofs,y=errH] {data/nontrivial/P1.txt}
node[pos=0.7,pin={[pin distance=0.15cm]-90:{\errH}}] {};
\plot[thick,red  ,mark=x     ] table[x=dofs,y=errB] {data/nontrivial/P1.txt}
node[pos=0.1,pin={[pin distance=0.7cm]-90:{\errB}}] {};
\plot[thick,black,mark=square] table[x=dofs,y=est ] {data/nontrivial/P1.txt}
node[pos=0.9,pin={[pin distance=0.5cm] 90:{\est }}] {};

\plot [dashed,domain=3e4:5e5] {16*x^(-2./3.)};

\SlopeTriangle{0.55}{-0.2}{0.1}{-0.66666}{$\ndofs^{-2/3}$}{}

\end{axis}
\end{tikzpicture}
\end{minipage}
\begin{minipage}{.45\linewidth}
\begin{tikzpicture}
\begin{axis}
[
width=\textwidth,
xmode = log,
ymin = 0.8,
ymax = 1.2,
xlabel={\ndof}
]

\plot[thick,black,mark=square] table[x=dofs,y expr=\thisrow{est}/\thisrow{errH}]
{data/nontrivial/P1.txt}
node[pos=0.5,pin={[pin distance=0.5cm]-90:{\est/\errH}}] {};

\end{axis}
\end{tikzpicture}
\end{minipage}
\caption{Non-trivial cohomology example, $p=1$}
\label{figure_nontrivial_P1}
\end{figure}

Figures \ref{figure_nontrivial_P0} and \ref{figure_nontrivial_P1}
show the behaviour of the post-processed induction field and the
estimator for $p=0$ and $1$ respectively. The behaviour is similar
than in the trivial cohomology example, with the difference that the
error is slightly underestimated on the first mesh. This is to be
expected since the data oscillation term is not included and starting
mesh is coarser in this example than in the trivial topology case.

\begin{figure}
\begin{minipage}{.45\linewidth}
\begin{tikzpicture}
\begin{axis}
[
width=\textwidth,
ymode = log,
xlabel={\ndof},
xtick={10,20,30,40,50},
xticklabels={1,8,27,64,125},
scaled x ticks={real:1.0},
xtick scale label code/.code={$\cdot 10^3$}]
]

\plot[thick,blue ,mark=o     ] table[x expr=\thisrow{dofs}^(1./3.),y=errH] {data/nontrivial/H.txt}
node[pos=0.4,pin={[pin distance=0.7cm]-90:{\errH}}] {};
\plot[thick,red  ,mark=x     ] table[x expr=\thisrow{dofs}^(1./3.),y=errB] {data/nontrivial/H.txt}
node[pos=0.1,pin={[pin distance=0.7cm]-90:{\errB}}] {};
\plot[thick,black,mark=square] table[x expr=\thisrow{dofs}^(1./3.),y=est ] {data/nontrivial/H.txt}
node[pos=0.9,pin={[pin distance=0.5cm] 90:{\est }}] {};

\plot [dashed,domain=10:50] {5*exp(-0.25*x)};

\SlopeTriangle{0.55}{-0.2}{0.15}{-0.25}{$e^{-0.25\ndofs^{1/3}}$}{}

\end{axis}
\end{tikzpicture}
\end{minipage}
\begin{minipage}{.45\linewidth}
\begin{tikzpicture}
\begin{axis}
[
width=\textwidth,
ymin = 0.9,
ymax = 1.1,
xlabel={\ndof},
xtick={10,20,30,40,50},
xticklabels={1,8,27,64,125},
scaled x ticks={real:1.0},
xtick scale label code/.code={$\cdot 10^3$}]
]

\plot[thick,black,mark=square]
table[x expr=\thisrow{dofs}^(1./3.),y expr=\thisrow{est}/\thisrow{errH}]
{data/nontrivial/H.txt}
node[pos=0.5,pin={[pin distance=0.5cm]90:{\est/\errH}}] {};

\end{axis}
\end{tikzpicture}
\end{minipage}
\caption{Non-trivial cohomology example, $p$-convergence}
\label{figure_nontrivial_H}
\end{figure}
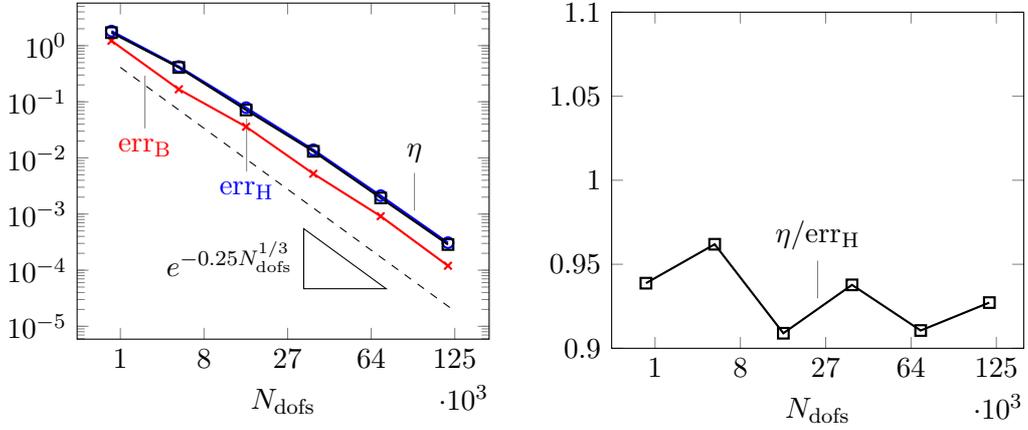

Figure \ref{figure_nontrivial_H} displays a $p$-convergence result
on a fixed mesh with 1043 tetrahedra ($h_{\max} = 1.19$). As before,
we observe that the estimator is $p$-robust. The error is slightly
underestimated due to the coarseness of the mesh.

\subsection{Adaptivity}

We close this section with an example where the estimator is employed
to drive an adaptive mesh refinement process with D\"orfler marking
\cite{dorfler_1996a}. After the solution and the estimator have been
computed, the elements are sorted according to the estimator values,
and we select a minimal set $\CM_h$ of elements such that
\begin{equation*}
\sum_{K \in \CM_h} \eta_K^2 \leq \theta^2 \sum_{K \in \CT_h} \eta_K^2,
\qquad
\theta \eq 0.05.
\end{equation*}
The software {\tt mmg3d} \cite{dobrzynski_2012a} is used to iteratively refine
the mesh locally around the marked elements $K \in \CM_h$, with the procedure
detailed in \cite[Section 6.1.4]{chaumontfrelet_2023a}.

The computational domain is the ``L-brick'' $\Omega \eq L \times (0,2)$
where $L = (-2,2)^2 \setminus (0,2) \times (-2,0)$. The right-hand side is given by
\begin{equation*}
\BJ
\eq
\left (
\begin{array}{c}
 (\bx_2-\by_2)
\\
-(\bx_1-\by_1)
\\
0
\end{array}
\right )
e^{-|\bx-\by|^2/\sigma^2},
\end{equation*}
with $\by = (-1,1,1)$ and $\sigma = 0.4$. Notice that this right-hand side
is divergence-free but does exactly satisfy the boundary condition
$\BJ \cdot \bn = 0$ on $\partial \Omega$. However, the normal trace is so
small that the discrepancy does not affect the numerical results for the mesh
sizes considered. The analytical solution $\BH$ is not available here, so
that a reference solution $\widetilde \BH$ is computed for each mesh by increasing
the polynomial degree by $1$.

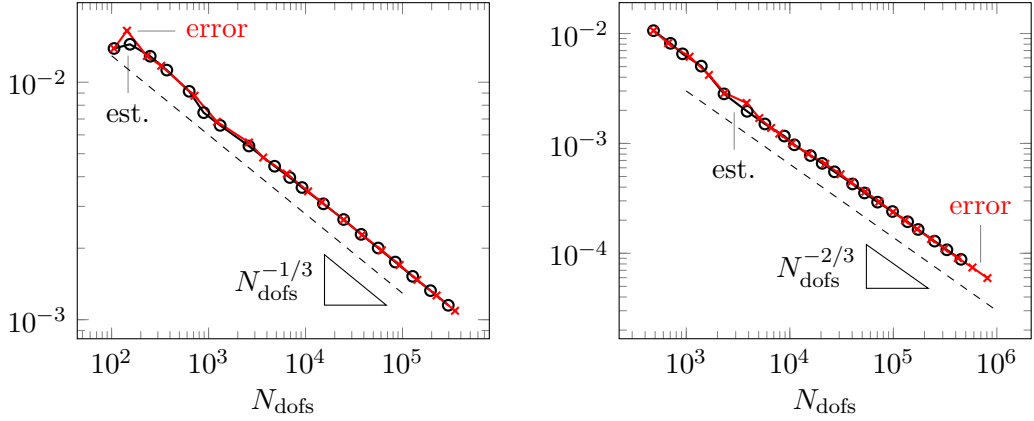
\begin{figure}
\begin{minipage}{.45\linewidth}
\begin{tikzpicture}
\begin{axis}
[
width=\textwidth,
ymode = log,
xmode = log,
xlabel = {\ndof}
]

\plot[black,thick,mark=o] table[x=dofs,y=errH] {data/adaptivity/P0.txt}
node[pos=0.04,pin=-90:{est.}] {};
\plot[red  ,thick,mark=x] table[x=dofs,y=errH] {data/adaptivity/P0_cheat.txt}
node[pos=0.04,pin=0:{error}] {};

\plot [dashed,domain=1e2:1e5] {0.06*x^(-1./3.)};

\SlopeTriangle{0.6}{-0.15}{0.1}{-0.333333}{$\ndofs^{-1/3}$}{}

\end{axis}
\end{tikzpicture}
\end{minipage}
\begin{minipage}{.45\linewidth}
\begin{tikzpicture}
\begin{axis}
[
width=\textwidth,
ymode = log,
xmode = log,
xlabel = {\ndof}
]

\plot[black,thick,mark=o] table[x=dofs,y=errH] {data/adaptivity/P1.txt}
node[pos=0.28,pin=-90:{est.}] {};
\plot[red  ,thick,mark=x] table[x=dofs,y=errH] {data/adaptivity/P1_cheat.txt}
node[pos=0.98,pin=90:{error}] {};

\plot [dashed,domain=1e3:1e6] {0.3*x^(-2./3.)};

\SlopeTriangle{0.6}{-0.15}{0.15}{-0.66666}{$\ndofs^{-2/3}$}{}

\end{axis}
\end{tikzpicture}
\end{minipage}
\caption{$\merrH$ in error-driven and estimator-driven adaptivity}
\label{figure_adaptivity_comparison}
\end{figure}

The adaptive refinement process is started with an initial mesh of 56 tetrahedra
($h_{\max} = 2.83$). We perform two simulations with $p=0$ and $p=1$. In Figure
\ref{figure_adaptivity_comparison}, we compare the convergence rates obtained
using the estimator $\eta_K$ and the ``true'' error $\|\widetilde \BH-\BH_h\|_K$
to mark the elements in the adaptivty refinements. The results are extremely similar,
which shows that the estimator is perfectly suited to drive the refinements. The optimal
convergence rate of $N_{\rm dofs}^{p/3}$ are also achieved.

Figures \ref{figure_adaptivity_P0} and \ref{figure_adaptivity_P1} display
the behaviour of the error, the estimator and the post-processed numerical
solution throughout the adaptive process. As before, we observe that the
post-processed solution is more accurate. The effectivity index is very close
to one, except on the coarsest mesh where the error is underestimated due to
the absence of the data oscillation term.

Figures \ref{figure_adaptivity_iterations_P0_0005}, \ref{figure_adaptivity_iterations_P0_0015},
\ref{figure_adaptivity_iterations_P1_0010} and \ref{figure_adaptivity_iterations_P1_0020}
represent the meshes obtained at different stages of the iterative refinement process.
Specifically the top size of the L-brick is represented on the top of the figures, whereas
the two faces sharing the reentrant edges are represented at the bottom.
The agreement between the estimator $\eta_K$ and the local error $\|\widetilde \BH-\BH_h\|_K$
is excellent. We also observe that the mesh is locally refined as expected: close to
the source center $\by$ and along the reentrant edge.

\begin{figure}
\begin{minipage}{.45\linewidth}
\begin{tikzpicture}
\begin{axis}
[
width=\textwidth,
ymode = log,
xmode = log,
xlabel = {\ndof}
]

\plot[thick,blue ,mark=o     ] table[x=dofs,y=errH] {data/adaptivity/P0.txt}
node[pos=0.9,pin={[pin distance=0.1cm]-90:{\errH}}] {};
\plot[thick,red  ,mark=x     ] table[x=dofs,y=errB] {data/adaptivity/P0.txt}
node[pos=0.3,pin={[pin distance=0.6cm]-90:{\errB}}] {};
\plot[thick,black,mark=square] table[x=dofs,y=est ] {data/adaptivity/P0.txt}
node[pos=0.9,pin={[pin distance=0.5cm] 90:{\est }}] {};

\plot [dashed,domain=1e2:2e5] {0.02*x^(-1./3.)};

\SlopeTriangle{0.6}{-0.2}{0.1}{-0.333333}{$\ndofs^{-1/3}$}{}

\end{axis}
\end{tikzpicture}
\end{minipage}
\begin{minipage}{.45\linewidth}
\begin{tikzpicture}
\begin{axis}
[
width=\textwidth,
xmode = log,
ymin = 0.1,
ymax = 1.9,
xlabel = {\ndof}
]

\plot[thick,black,mark=square] table[x=dofs,y expr=\thisrow{est}/\thisrow{errH}]
{data/adaptivity/P0.txt}
node[pos=0.5,pin={[pin distance=0.5cm]-90:{\est/\errH}}] {};

\end{axis}
\end{tikzpicture}
\end{minipage}
\caption{Adaptivity example, $p=0$}
\label{figure_adaptivity_P0}
\end{figure}

\begin{figure}
\begin{minipage}{.45\linewidth}
\begin{tikzpicture}
\begin{axis}
[
width=\textwidth,
ymode = log,
xmode = log,
xlabel = {\ndof}
]

\plot[thick,blue ,mark=o     ] table[x=dofs,y=errH] {data/adaptivity/P1.txt}
node[pos=0.0,pin={[pin distance=0.5cm]0:{\errH}}] {};
\plot[thick,red  ,mark=x     ] table[x=dofs,y=errB] {data/adaptivity/P1.txt}
node[pos=0.3,pin={[pin distance=0.5cm]-90:{\errB}}] {};
\plot[thick,black,mark=square] table[x=dofs,y=est ] {data/adaptivity/P1.txt}
node[pos=0.9,pin={[pin distance=0.5cm] 90:{\est }}] {};

\plot [dashed,domain=1e3:2e5] {0.1*x^(-2./3.)};

\SlopeTriangle{0.6}{-0.2}{0.1}{-0.66666}{$\ndofs^{-2/3}$}{}

\end{axis}
\end{tikzpicture}
\end{minipage}
\begin{minipage}{.45\linewidth}
\begin{tikzpicture}
\begin{axis}
[
width=\textwidth,
xmode = log,
ymin = 0.5,
ymax = 1.5,
xlabel = {\ndof}
]

\plot[thick,black,mark=square] table[x=dofs,y expr=\thisrow{est}/\thisrow{errH}]
{data/adaptivity/P1.txt}
node[pos=0.5,pin={[pin distance=0.5cm]-90:{\est/\errH}}] {};

\end{axis}
\end{tikzpicture}
\end{minipage}
\caption{Adaptivity example, $p=1$}
\label{figure_adaptivity_P1}
\end{figure}

\input{figures/adaptivity/iterations/P0_0005}
\input{figures/adaptivity/iterations/P0_0015}
\input{figures/adaptivity/iterations/P1_0010}
\input{figures/adaptivity/iterations/P1_0020}

\section{Conclusion}
\label{section_conclusion}

We propose an a posteriori error estimator for mixed finite element discretizations
of the curl-curl problem. Following the framework of Prager--Synge estimates, we
engineer a potential reconstruction leading to guaranteed error bounds. A key novelty
of our approach is that instead of explicitly constructing the potential, we actually
build its curl through patch-wise divergence-constrained minimization problems. This
construction leads to an estimator that is locally efficient and polynomial-degree-robust.
It is noteworthy that the reconstruction technique proposed here for the
curl-curl problem is similar to the construction of equilibrated fluxes for
the Poisson problem.

We highlight the properties of the reconstructed field and associated estimator with
selected numerical examples. These examples show that the reconstructed induction field
is typically more accurate than the discrete magnetic field output by the finite element
scheme, although not being super-convergent. The key theoretical properties are also
observed numerically: the proposed estimator is efficient, polynomial-degree-robust,
and fully reliable up to data oscillations. We also employ the estimator to drive
adaptive mesh refinements in a domain with a reentrant edge. This example indicates
that the estimator is suited for adaptivity purposes: the meshes produced are correctly
refined, leading to optimal convergence rates.

\bibliographystyle{amsplain}
\bibliography{bibliography}

\end{document}